\documentclass{article}
\usepackage{amsmath,amssymb,amsthm}
\usepackage[headings]{fullpage}

\usepackage{subfig,xcolor,graphicx}
\usepackage{algorithm,algorithmic}

\ifx \eqref \undefined
\newcommand{\eqref}[1]{(\ref{#1})}
\fi

\newtheorem{thrm}{Theorem}
\newtheorem{lem}[thrm]{Lemma}
\newtheorem{cor}[thrm]{Corollary}

\newtheorem{conj}[thrm]{Conjecture}

\newcommand{\bvec}[1] {{\mathbf {#1}}}

\newcommand{\BPC}{{P_{\text{cr}}^{\text{\tiny B}}}}

\newcommand{\NR}{{R}}
\newcommand{\BR}{{R_{\text{\tiny B}}}}

\newcommand{\BRone}{{R_{\text{\tiny B,1D}}}}

\newcommand{\NS}{{S}}
\newcommand{\BS}{{S_{\text{\tiny B}}}}

\newcommand{\NQ}{{Q}}
\newcommand{\BQ}{{Q_{\text{\tiny B}}}}
\newcommand{\alphan}{{\alpha_{\text{\tiny NLS}}}}
\newcommand{\alphab}{{\alpha_{\text{\tiny B}}}}

\newcommand{\TCrit}{{T_{\rm c}}}

\newcommand{\rmax}{ { {\it r_{\max}}  } }

\newcommand{\Schrodinger}{{Schr\"odinger }}

\ifx \DeclareMathOperator \undefined

\else

\fi

\newcommand{\mycaption}[1]{\parbox{0.95\textwidth}{\caption{{#1}}}}

\newcommand{\Real}{{\mathbb R}}
\newcommand{\abs}[1]{{\left\vert{#1}\right\vert}}
\newcommand{\norm}[1]{{\left\Vert{#1}\right\Vert}}

\newcommand{\LN}{{ \mathit l }}

\newcommand{\myfw}{ 0.4 }
\newcommand{\myfwb}{ 0.35 }

\begin{document}

\renewcommand{\thesubfigure}{{\Alph{subfigure}}}

\title{Ring-type singular solutions of the biharmonic nonlinear \Schrodinger
equation}

\makeatletter
\author{{G. Baruch}, {G. Fibich}$^*$, {E. Mandelbaum},\\
    School of Mathematical Sciences, Tel Aviv University, Tel Aviv 69978,
    Israel\\
    $^*$Corresponding author, fibich@math.tau.ac.il
}
\makeatother
\maketitle

\begin{abstract}
    We present new singular solutions of the biharmonic nonlinear
    \Schrodinger equation \[
        i\psi_t(t,\bvec{x}) - \Delta^2\psi + |\psi|^{2\sigma}\psi =0,
        \qquad  \bvec{x}\in \mathbb{R}^d,
        \qquad  4/d\le\sigma\le4.
    \]
    These solutions collapse with the quasi self-similar ring
    profile~$\psi_{\BQ}$, where \[
        \abs{\psi_{\BQ}(t,r)}
            \sim \frac{1}{L^{2/\sigma}(t)}
            \BQ\left( \frac {r-\rmax(t)}{L(t)} \right),
            \qquad r = \abs{\bvec{x}},
    \]
    $L(t)$ is the ring width that vanishes at singularity,~$
        \rmax(t)\sim r_0L^{\alpha}(t)
    $ is the ring radius, and~$
        \alpha
            = \frac{4-\sigma}{\sigma(d-1)}
    $
    The blowup rate of these solutions is~$ \frac{1}{3+\alpha} $
    for~$4/d\le\sigma<4$, and slightly faster than~$1/4$ for~$\sigma=4$.
    These solutions are analogous to the ring-type solutions of the
    nonlinear \Schrodinger equation.
\end{abstract}

%\tableofcontents
%\listoffigures

\section{\label{sec:intro}Introduction}

%\subsection{Singular solutions of the biharmonic nonlinear Schr\"odinger
%equation}
The focusing nonlinear \Schrodinger equation (NLS)
\begin{equation} \label{eq:NLS}
    i\psi_ t( t,{\bf x})+\Delta\psi+|\psi|^{2\sigma}\psi=0,\qquad
    \psi(0,{\bf x})=\psi_0({\bf x})
        \in H^1( \mathbb{R}^d ),
\end{equation}
where~${\bf x}\in\mathbb{R}^d$ and~$
    \Delta=\partial_{x_1x_1}+\cdots+\partial_{x_dx_d},
$ admits solutions that become singular at a finite time, i.e.,~$
    \lim_{t\to\TCrit}\norm{\psi}_{H^1} = \infty,
$
where~$0\le t\le\TCrit$.
Until a few years ago, all known singular solutions of the NLS were peak-type.
By this, we mean that if assume radial symmetry, and denote the
location of maximal amplitude by
\[
    \rmax(t) = \arg \displaystyle\max_r\abs{\psi},
    \qquad r=\sqrt{x_1^2+\dots+x_d^2},
\]
then~$\rmax(t)\equiv0$ for~$0\le t \le \TCrit$, i.e., the solution peak is
attained at~$r=0$.
In recent years, however, new singular solutions of the NLS were found, which
are ring-type, i.e.,~$\rmax(t)>0$ for~$0\le t < \TCrit$.

In this study, we consider the focusing {\em biharmonic} nonlinear
\Schrodinger equation (BNLS)
\begin{equation}    \label{eq:BNLS}
    i\psi_t(t,\bvec{x}) - \Delta^2\psi + \left|\psi\right|^{2\sigma}\psi = 0,
    \qquad \psi(0,\bvec{x}) = \psi_0(\bvec{x})\in H^2(\Real^d),
\end{equation}
where~$\Delta^2$ is the biharmonic operator.
Singular peak-type solutions of the BNLS have been studied
in~\cite{Fibich_Ilan_George_BNLS:2002,Baruch_Fibich_Mandelbaum:2009a,Baruch_Fibich:2010a}.
Singular ring-type solutions of the BNLS with~$\sigma>4$ were studied
in~\cite{Baruch_Fibich_Gavish:2009}.
The goal of this work is to find and characterize singular ring-type solutions
of the BNLS with~$4/d\le\sigma\le4$.

\subsection{\label{ssec:NLS}Singular solutions of the nonlinear \Schrodinger equation (NLS) - review}

The NLS~\eqref{eq:NLS} is called 
{\em subcritical} if $\sigma d<2$.
In this case, all solutions exist globally.
In contrast, solutions of the {\em critical}~($\sigma d=2$) and
{\em supercritical}~($\sigma d>2$) NLS can become singular at a finite time.

Until a few years ago, the only known singular NLS solutions were peak-type.
In the critical case~$\sigma d=2$, it has been rigorously shown~\cite{Merle-03}
that peak-type solutions are self-similar near the singularity,
i.e.,~$\psi\sim\psi_\NR$, where
\begin{equation*}
    \psi_\NR(t,r)=
        \frac{1}{L^{d/2}(t)}
        \NR\left( \frac r{L(t)} \right) e^{i\int_0^t\frac{ds}{L^2(s)}},
\end{equation*}
and~$r=\abs{\bvec{x}}$.
The self-similar profile~$\NR(\rho)$ is the ground state of the standing-wave
equation\[
    \NR^{\prime\prime}(\rho) +\frac{d-1}{\rho}\NR^\prime -\NR+|\NR|^{4/d}\NR=0.
\]
Since~$\NR$ attains its global maximum at~$\rho=0$, $\psi_\NR$ is a peak-type
profile.
The blowup rate of~$L(t)$ is given by the {\em loglog law}
\begin{equation}
    L(t) \sim \left(
        \frac{2\pi(\TCrit-t)}{\log\log1/(\TCrit-t)}
    \right)^{1/2},
    \qquad t\to \TCrit.
\label{eq:logloglaw}
\end{equation}

In the supercritical case~$\sigma d>2$, the rigorous theory is far less
developed.
However, formal calculations and numerical simulations~\cite{Sulem-99} suggest
that peak-type solutions of the supercritical NLS collapse with the
self-similar~$\psi_\NS$ profile, i.e.,~$\psi\sim\psi_\NS$, where
\begin{subequations}    \label{eq:intro_psiS}
    \begin{equation}    \label{eq:intro_psiS_profile}
        \psi_\NS(t,r) = \frac{1}{L^{1/\sigma}(t)}
            \NS\left(\rho\right)e^{i\tau},
    \end{equation}
    \begin{equation}    \label{eq:intro_psiS_profile2}
        \tau=\int_0^t\frac{ds}{L^2(s)},
        \qquad  \rho=\frac{r}{L(t)},
    \end{equation}
    and~$\NS(\rho)$ is the zero-Hamiltonian, monotonically-decreasing solution of
    the nonlinear eigenvalue problem
    \begin{equation} \label{eq:intro_psiS_ODE4S}
            S^{\prime\prime}(\rho) +\frac{d-1}{\rho}S^\prime - S 
            +i\frac {\kappa^2} 2 \left( \frac 1\sigma S + \rho S^\prime \right)
            + |S|^{2\sigma}S = 0,
            \qquad S^\prime(0)=0,
    \end{equation}
    where~$\kappa$ is the eigenvalue.
    Since~$|\NS(\rho)|$ attains its global maximum at~$\rho=0$, $\psi_\NS$ is a
    peak-type profile.
    The blowup rate of~$L(t)$ is a square-root, i.e.,
    \begin{equation}    \label{eq:intro_psiS_blowuprate}
        L(t) \sim \kappa\sqrt{\TCrit-t},
        \qquad t\to \TCrit,
    \end{equation}
    where~$\kappa>0$ is the eigenvalue of~\eqref{eq:intro_psiS_ODE4S}.
\end{subequations}

In the last few years, new singular solutions of the NLS were discovered,
which are
ring-type~\cite{Gprofile-05,SC_rings-07,Raphael-06,Raphael-08,Vortex_rings-08}.
In particular, in~\cite{SC_rings-07}, Fibich, Gavish and Wang showed that the
NLS~\eqref{eq:NLS} with~$d>1$ and~$\frac2d\le\sigma\le2$ admits singular
ring-type solutions that collapse with the~$\psi_{\NQ}$ profile,
i.e.,~$\psi\sim\psi_{\NQ}$, where
\begin{subequations}    \label{eq:psi_QN}
    \begin{equation}    \label{eq:psi_QN_profile}
        \psi_{\NQ}(t,r) =
            \frac{1}{L^{1/\sigma}(t)} \NQ(\rho)
            e^{
                i\tau
                +i\alphan\frac{L_t}{4L}r^2
                +i(1-\alphan)\frac{L_t}{4L}(r-\rmax(t))^2
            },
    \end{equation}
    \begin{equation}
        \tau = \int_0^t\frac{ds}{L^2(s)},
        \qquad  \rho=\frac{r-\rmax( t)}{L(t)},
        \qquad  \rmax(t)\sim r_0 L^\alphan(t),
    \end{equation}
    and
    \begin{equation}    \label{eq:psi_QN_alphan}
        \alphan = \frac{2-\sigma}{\sigma(d-1)}
        = 1 - \frac{\sigma d-2}{\sigma(d-1)}.
    \end{equation}
\end{subequations}
The self-similar profile~$\NQ$ attains its global maximum at~$\rho=0$.
Hence,~$\rmax(t)$ is the ring radius and~$L(t)$ is the ring width, see
Figure~\ref{fig:rm_def}.

A unique feature of the~$\psi_{\NQ}$ profile~\eqref{eq:psi_QN_profile} is the
linear combination of the two radial phase terms.
The first phase term~$
    \alphan\frac{L_t}{4L}r^2
$ describes focusing towards~$r=0$, and is the manifestation of the shrinking
of the ring radius~$\rmax$ to zero.
The second term~$
    (1-\alphan)\frac{L_t}{4L}(r-\rmax(t))^2
$ describes focusing towards~$r=\rmax$, and is the manifestation of the
shrinking of the ring width~$L(t)$ to zero.
The discovery of this ``double-lens'' ansatz was the key stage in the
asymptotic analysis of the~$\psi_{\NQ}$ profile, which enabled the calculation of
the shrinking rate~$\alphan$, see~\eqref{eq:psi_QN_alphan} and the blowup 
rate~$p$, see~\eqref{eq:NLS_shrinking_p}.

\begin{figure}\begin{center}
    \scalebox{.4}{\input{ring_illustration.pstex_t}}
    \mycaption{Illustration of ring radius~$\rmax(t)$ and width~$L( t)$.
    \label{fig:rm_def}}
\end{center}\end{figure}

\begin{figure}[h]
    \begin{center}
        \scalebox{0.7}{\input{NLS_ring_class.pstex_t}}

        \mycaption{ \label{fig:NLS_class_diagram}
            Classification of singular ring-type solutions of the NLS, as a
            function of~$\sigma$ and~$d$.
            A)~subcritical case - no singular solutions exist.
            B)~$\sigma d=2$: equal-rate~$\psi_{\NQ}$
            solutions~\cite{Gprofile-05}.
            C)~$2/d<\sigma<2$: shrinking~$\psi_{\NQ}$
            solutions~\cite{SC_rings-07}.
            D)~$\sigma=2$: standing~$\psi_{\NQ}$
            solutions~\cite{SC_rings-07,Raphael-06,Raphael-08}.
            E)~$\sigma>2$: standing non-$\psi_{\NQ}$
            rings~\cite{Baruch_Fibich_Gavish:2009}.
        }
    \end{center}
\end{figure}
The NLS ring-type singular solutions can be classified as follows, see
Figure~\ref{fig:NLS_class_diagram}:
\begin{enumerate}
    \renewcommand{\theenumi}{\Alph{enumi}}
    \item In the subcritical case~($\sigma d<2$), all NLS solutions exist
        globally, hence no singular ring-type solutions exist.
    \item The critical case~$\sigma d=2$ corresponds to~$\alphan=1$.
        Since~$\rmax(t)\sim r_0 L( t)$, these solutions undergo an
        {\em equal-rate collapse}, i.e., the ring radius goes to zero at the
        same rate as~$L(t)$.
        The blowup rate of~$L(t)$ is a square root.
    \item The supercritical case~$2/d<\sigma<2$ corresponds to $0<\alphan<1$.
        Therefore, the ring radius~$\rmax(t)\sim r_0 L^\alphan(t)$ decays to
        zero, but at a slower rate than~$L(t)$.
        The blowup rate of~$L(t)$ is
        \begin{equation}\label{eq:blowup_rate_with_p}
            L(t)\sim \kappa(\TCrit-t)^p,
        \end{equation}
        where
        \begin{equation}    \label{eq:NLS_shrinking_p}
            p = \frac{1}{1+\alphan}
            = \frac{1}{2-\frac{\sigma d-2}{\sigma (d-1)}}
        \end{equation}
    \item The supercritical case~$\sigma=2$ corresponds to~$\alphan=0$, i.e.,~$
            \displaystyle \lim_{t\to\TCrit}\rmax(t) = \rmax(\TCrit)>0
        $.
        Therefore, the solution becomes singular on the d-dimensional
        sphere~$|\bvec{x}|=\rmax(\TCrit)$, rather than at a point.
        The blowup profile~$\psi_{\NQ}$ is equal to that of peak-type solutions of
        the $1D$ critical NLS, and the blowup rate is given by the loglog
        law~\eqref{eq:logloglaw}.
    \item The case~$\sigma>2$ also corresponds to a standing ring.
        The asymptotic profile is not given by~$\psi_{\NQ}$, however, but rather by
        the asymptotic profile of peak-type solutions of the $1D$ supercritical
        NLS.
        The blowup rate is a square root~\cite{Baruch_Fibich_Gavish:2009}.
\end{enumerate}
\noindent
Thus, NLS ring-type singular solutions are {\em shrinking}
(i.e.,~$\lim_{t\to \TCrit}{\rmax(t)}=0$) for~$\frac2d\le \sigma <2$
(cases B and C), and {\em standing}
(i.e.,~$0<\lim_{t\to \TCrit}{\rmax(t)}<\infty$) for~$\sigma\ge2$
(cases D and E).

\subsection{\label{ssec:intro_BNLS}Singular solutions of the biharmonic NLS -
review}

The BNLS equation~\eqref{eq:BNLS} is called
{\em subcritical} if $\sigma d<4$,
{\em supercritical} if $\sigma d>4$, and
{\em critical} if $\sigma d=4$.
In the critical case, equation~\eqref{eq:BNLS} can be rewritten as
\begin{equation}    \label{eq:CBNLS}
    i\psi_t(t,\bvec{x}) - \Delta^2\psi + \left|\psi\right|^{8/d}\psi = 0,
    \qquad \psi(0,\bvec{x}) = \psi_0(\bvec{x})\in H^2(\Real^d).
\end{equation}

The BNLS conserves the ``power'' ($L^2$ norm), i.e., \[
	P(t) \equiv P(0),\qquad
	P(t) = \norm{\psi(t)}_2^2,
\]
and the Hamiltonian
\begin{equation}	\label{eq:Hamiltonian_conservation}
	H(t) \equiv H(0),\qquad 
	H[\psi(t)]  =  
		\norm{\Delta\psi}_2^2 
			-{\textstyle\frac{1}{1+\sigma}}
            \norm{\psi}_{2(\sigma+1)}^{2(\sigma+1)}.
\end{equation}

In the radially-symmetric case, the BNLS equation~\eqref{eq:BNLS} reduces
to
\begin{equation}    \label{eq:radial_BNLS}
    i\psi_t(t,r) - \Delta^2_r\psi + \left|\psi\right|^{2\sigma}\psi = 0,
    \qquad \psi(0,r) = \psi_0(r),
\end{equation}
where
\begin{equation}    \label{eq:radial_bi_Laplacian}
    \Delta_r^2 =
        \partial_r^4
        +\frac{2(d-1)}{r}\partial_r^3
        +\frac{(d-1)(d-3)}{r^2}\partial_r^2
        -\frac{(d-1)(d-3)}{r^3}\partial_r
\end{equation}
is the radial biharmonic operator.
In particular, the radially-symmetric critical BNLS is given by
\begin{equation}    \label{eq:radial_CBNLS}
    i\psi_t(t,r) - \Delta^2_r\psi + \left|\psi\right|^{8/d}\psi = 0,
    \qquad \psi(0,r) = \psi_0(r).
\end{equation}

All solutions of the subcritical BNLS exist
globally in~$H^2$~\cite{Fibich_Ilan_George_BNLS:2002}.
In the critical case, they exist globally if the input power is below the
critical power:
\begin{thrm}[\cite{Fibich_Ilan_George_BNLS:2002}]    \label{thrm:GE_critical}
    Let
    $ \norm{\psi_0}_2^2 < \BPC$, where~$\BPC = \norm{\BR}_2^2$,
    and~$\BR(\rho)$ is the ground state of the standing wave equation
    \begin{equation}    \label{eq:standing_wave}
        -\Delta_{\rho}^2\BR(\rho)-\BR+\abs{\BR}^{2\sigma}\BR = 0,
    \end{equation}
    with~$\sigma=4/d$.
    Then, the solution of the critical focusing BNLS~\eqref{eq:CBNLS} exists
    globally.
\end{thrm} % thrm:GE_critical
\noindent
Numerical simulations~\cite{Fibich_Ilan_George_BNLS:2002,%
Baruch_Fibich_Mandelbaum:2009a} indicate that solutions of the critical and
supercritical BNLS can become singular at a finite time, i.e.,~$
    \lim_{ t\to \TCrit}\|\psi\|_{H^2}=\infty,
$ where $0< \TCrit<\infty$.
At present, however, there is no rigorous proof that the BNLS admits singular
solutions, whether peak-type or ring-type.

In~\cite{Baruch_Fibich_Mandelbaum:2009a}, we rigorously proved that the blowup
rate of all $H^2$ singular solutions of the critical BNLS is bounded by a
quartic root:
\begin{thrm}    \label{thrm:low-bound}
Let~$\psi$ be a solution of the critical BNLS~(\ref{eq:CBNLS}) that
becomes singular at~$t=\TCrit<\infty$, and let~$
    \LN(t) = \norm{\Delta\psi}_2^{-1/2}
$.
Then,~$\exists K=K(\norm{\psi_0}_2)>0$ such that  \[
    \LN(t) \leq  K(\TCrit-t)^{1/4},
    \qquad 0\leq t<\TCrit.
\]
\end{thrm}

\noindent
We also proved that all singular solutions are quasi self-similar:
\begin{thrm}    \label{thrm:self-similarity}
Let~$d\geq 2$, and let~$\psi(t,r)$ be a solution of the radially-symmetric
critical BNLS~(\ref{eq:radial_CBNLS}) with initial conditions~$
    \psi_0(r)\in H^2_{\rm radial}
$, that becomes singular at~$t=\TCrit<\infty$.
Let~$ \LN(t) = \norm{\Delta\psi}_2^{-1/2}$, and let \[
    S(\psi)(t,r)  =  \LN^{d/2}(t)\psi(t,\, \LN(t)r).
\]
Then, for any sequence~$t'_k\to \TCrit$, there is a subsequence~$t_k$, such
that $S(\psi)(t_k,r)\to\Psi(r)$ strongly in~$L^q$, for all~$q$ such that%
%\footnote{
%   In fact,~$q=2(\sigma+1)$, where~$\sigma$ is in the $H^2$-subcritical
%   regime~\eqref{eq:admissible-range}.
%}%
\begin{equation}    \label{eq:thrm-self-q}
    %\begin{cases}
    \left\{ \begin{array}{cc}
        2<q<\infty  &   2\leq d\leq 4, \\
        2<q<\frac{2d}{d-4}  &   4<d .
    \end{array} \right.
    %\end{cases}
\end{equation}%
%In addition,~$\norm{\Psi}_2^2\geq \norm{R}_2^2$, where~$R$ is the ground
%state of equation~\eqref{eq:stationary_state}.
\end{thrm}

Since the $L^2$-norm of~$S(\psi)$ is conserved, and the convergence
of~$S(\psi)$ to~$\Psi$ is in~$L^q$ with~$q>2$, the solution becomes
self-similar in the singular region (the collapsing core), but not everywhere.
Consequentially, the solution has the power-concentration property, whereby a
finite amount of power enters the singularity point, i.e., \[
    \lim_{\varepsilon\to0+}
    \liminf_{t\to\TCrit}
    \norm{\psi}_{L^2(r<\varepsilon)}^2 \ge \BPC,
\]
where~$\BPC$ is the critical power for
collapse~\cite{Baruch_Fibich_Mandelbaum:2009a,ChaeHongLee:2009}.

Peak-type singular solutions of the critical BNLS~\eqref{eq:radial_CBNLS} were
studied asymptotically and numerically
in~\cite{Fibich_Ilan_George_BNLS:2002,Baruch_Fibich_Mandelbaum:2009a}.
The asymptotic profile of these solutions is
\begin{equation}    \label{eq:intro_peak_CBNLS}
    \psi_{\BR}(t,r)
    = \frac{1}{L^{d/2}(t)}
        \BR\left( \frac r{L(t)} \right)
        e^{i\int^t\frac{1}{L^4(s)}ds},
\end{equation}
where~$\BR(\rho)$ is the ground state of~\eqref{eq:standing_wave}.
The blowup rate of~$L(t)$ is slightly faster than a quartic root, i.e.,
\begin{equation}    \label{eq:rate_slightly_mad}
        \lim_{t\to\TCrit}
            \frac{L(t)}{(\TCrit-t)^p}
        =
        %\begin{cases}
        \left\{\begin{array}{cc}
            0\quad  & p=\frac14 \\
            \infty  & p>\frac14
        \end{array}\right.
        %\end{cases} .
\end{equation}

Specifically, in the one-dimensional case, the quasi self-similar profile is
\begin{subequations}    \label{eqs:phi-QSS}
    \begin{equation}    \label{eq:phi-QSS}
        \psi_{R_{B,1D}}(t,x) =
            \frac{1}{L^{1/2}(t)} \BRone
            \left( \frac x{L(t)} \right)
            e^{\int^{t}\frac{1}{L^{4}(s)}ds},
    \end{equation}
    and~$\BRone$ is the ground state of%
%   \footnote{
%       Numerical methods for the calculation of the ground state
%       of~\eqref{eq:standing-1D} are given
%       in~\cite{Baruch_Fibich_Mandelbaum:2009a,Fibich_Ilan_George_BNLS:2002}.
%   }
    \begin{equation}    \label{eq:standing-1D}
        -\BR^{\prime\prime\prime\prime}(\xi)-\BR+\abs{\BR}^8\BR = 0.
    \end{equation}
\end{subequations}

Peak-type solutions of the supercritical BNLS~\eqref{eq:radial_BNLS} were
studied asymptotically and numerically in~\cite{Baruch_Fibich:2010a}.
The asymptotic profile of these solutions is 
\begin{equation}    \label{eq:intro_peak_SCBNLS}
    \psi_{\BS}(t,r)
    = \frac{1}{L^{2/\sigma}(t)}
        \BS\left( \frac r{L(t)} \right)
        e^{i\int^t\frac{1}{L^4(s)}ds},
\end{equation}
where~$\BS(\rho)$ is the zero-Hamiltonian solution of a nonlinear eigenvalue
problem
\begin{equation}    \label{eq:supercrit_peak_ODE-0}
    %\begin{gathered}
    -\BS(\rho) + i\frac{\kappa^4}4 \left(
            \frac{2}{\sigma}\BS + \rho \BS^\prime 
        \right)
        - \Delta_\rho \BS + \abs{\BS}^{2\sigma}\BS = 0, %\\
%    \BS(0)=1,
    \qquad \BS^\prime(0)=\BS^{\prime\prime\prime}(0)=0,
%    \qquad H[\BS]=0,
    %\end{gathered}
\end{equation}
and $\kappa$ is the eigenvalue.
%\begin{equation}    \label{eq:B_Hamiltonian}
%    H[\BS] = \int_{\rho=0}^{\infty} \left(
%            \abs{ \Delta_\rho \BS }^2 
%            -\frac1{1+\sigma} \abs{ \BS }^{2+2\sigma} 
%        \right)
%        \rho^{d-1}d\rho
%\end{equation}
The blowup rate is exactly~$p=1/4$, i.e.,
\begin{equation}    \label{eq:blowup_rate_14}
    L(t) \sim \kappa\left( \TCrit-t \right)^{1/4},
\end{equation}
where~$\kappa>0$ is the nonlinear eigenvalue
of~\eqref{eq:supercrit_peak_ODE-0}.

Ring-type singular solutions of the supercritical BNLS~\eqref{eq:radial_BNLS}
with~$\sigma>4$ were studied asymptotically and numerically
in~\cite{Baruch_Fibich_Gavish:2009}.
These solutions are standing rings, i.e.,~$
    \displaystyle\lim_{t\to\TCrit} \rmax(t)>0
$.
The self-similar profile of these standing-ring solutions is
\begin{equation}    \label{eq:intro_ring_SCBNLS}
    \psi_B(t,r)
    = \psi_{S_{B,1D}}(t,x=r-\rmax(t))
    = \frac{1}{L^{2/\sigma}(t)}
    S_{\text{B,1D}}\left( \frac {r-\rmax(t)}{L(t)} \right)
        e^{i\int^t\frac{1}{L^4(s)}ds},
\end{equation}
where~$\psi_{S_{\text{B,1D}}}(t,x)$, see~\eqref{eq:intro_peak_SCBNLS}, is the
profile of the peak-type singular solution of the one-dimensional supercritical
BNLS with the same value of~$\sigma$.
The blowup rate is given by~\eqref{eq:blowup_rate_14}.

\subsection{\label{ssec:intro_analogy}Analogy of NLS and BNLS}

\begin{table}
    \centering
    \begin{tabular}{|c||c|c|}
        \hline
        &   NLS &   BNLS \\
        \hline
        \multicolumn{3}{|c|}{peak-type solutions: critical
        case~\cite{Baruch_Fibich_Mandelbaum:2009a}} \\
        \hline
        \hline
        & $\sigma d=2$    &   $\sigma d=4$ \\
        \hline
        asymptotic profile &
        $
            \frac{1}{L^{1/\sigma}(t)}
            \NR\left( \frac r{L(t)} \right)
            e^{i\int_0^t\frac{ds}{L^2(s)}}
        $ &
        $
            \frac{1}{L^{2/\sigma}(t)}
            \BR\left( \frac r{L(t)} \right)
            e^{i\int_0^t\frac{ds}{L^4(s)}}
        $ \\
        \hline
        %\cline{2-3}
        blowup rate & slightly faster than~$1/2$ & 
            slightly faster than~$1/4$ \\
        \hline
        power concentration & yes & yes \\
        \hline
        \multicolumn{3}{|c|}{peak-type solutions: supercritical
        case~\cite{Baruch_Fibich:2010a}} \\
        \hline
        \hline
        & $\sigma d>2$    &   $\sigma d>4$ \\
        \hline
        asymptotic profile &
        $
            \frac{1}{L^{1/\sigma}(t)}
            \NS\left( \frac r{L(t)} \right) e^{i\int_0^t\frac{ds}{L^2(s)}}
        $ &
        $
            \frac{1}{L^{2/\sigma}(t)}
            \BS\left( \frac r{L(t)} \right) e^{i\int_0^t\frac{ds}{L^4(s)}}
        $ \\
        \hline
        blowup rate & $=1/2$ & $=1/4$ \\
        \hline
        power concentration & no & no \\
        \hline
        \multicolumn{3}{|c|}{``supercritical'' standing-ring
        solutions~\cite{Baruch_Fibich_Gavish:2009}} \\
        \hline
        \hline
            & $\sigma>2$    &   $\sigma>4$ \\
        \hline
        asymptotic profile & 
            same as $1D$ peak &
            same as $1D$ peak \\
        \hline
        blowup rate & 
            same as $1D$ peak &
            same as $1D$ peak \\
        \hline
    \end{tabular}
    \mycaption{
        A comparison of the properties of singular solutions of the NLS and
        BNLS.
        These properties are analogous, ``up to the change~$2\to4$''.
    }
    \label{tab:analogy}
\end{table}

Table~\ref{tab:analogy} lists the major findings of the previous
works~\cite{Fibich_Ilan_George_BNLS:2002,Baruch_Fibich_Mandelbaum:2009a,Baruch_Fibich_Gavish:2009}
on singular solutions of the BNLS, side by side with their NLS counterparts.
In all cases, the results for the BNLS mirror those of the NLS exactly,
``up to the change of~$2\to4$''.

We note that current BNLS theory is still missing a key feature in NLS theory,
which is the BNLS analogue of the quadratic radial phase terms of the
asymptotic profiles.
Therefore, our asymptotic analysis of the BNLS singular solutions produces
weaker results than those of~\cite{SC_rings-07}.
Hence, in this work we ``fill in'' the missing results by relying on the above
analogy of the NLS and BNLS, up to the change~$2\to4$.

\subsection{\label{ssec:intro_summary}Summary of results}

In this study, we consider ring-type singular solutions of the
BNLS~\eqref{eq:radial_BNLS} with~$4/d\leq \sigma \leq 4$.
We show numerically that such solutions exist, and are of the
form~$\psi(t,r)\sim\psi_{\BQ}(t,r)$, where
\begin{subequations}
    \label{eq:psi_QB}
    \begin{equation}    \label{}
        \abs{\psi_{\BQ}(t,r)}=
        \frac{1}{L^{2/\sigma}(t)} \BQ(\rho),
    \end{equation}
    \begin{equation}
        %\tau(t) = \int^{t}\frac{1}{L^4(s)}ds,
        \rho=\frac{r-\rmax(t)}{L(t)},
        \qquad  \rmax(t) \sim r_0 L^{\alpha}(t),
    \end{equation}
    and
    \begin{equation}    \label{eq:psi_QN_alpha_BNLS}
        \alpha 
        = \alphab
        = \frac{4-\sigma}{\sigma(d-1)}
        = 1-\frac{\sigma d-4}{\sigma(d-1)}.
    \end{equation}
\end{subequations}
The~$\psi_{\BQ}$ profile is the BNLS analogue of the~$\psi_{\NQ}$ profile of the
NLS.
Unlike the~$\psi_{\NQ}$ profile, however, we do not know the expression for the
double-lens ansatz of~$\psi_\BQ$.

In Section~\ref{sec:standing}, we consider the case~$\sigma=4$.
In this case~$\alpha=0$, i.e., the solution is a singular standing ring.
Informal asymptotic analysis and numerical simulations show that the blowup
profile is the self-similar profile
\begin{equation}    \label{eq:intro_ring_CBNLS}
    \psi_{\BQ}(t,r;\sigma=4)
    = \frac{1}{L^{1/2}(t)}
        \BRone\left( \frac {r-\rmax(t)}{L(t)} \right)
        e^{i\int^t\frac{1}{L^4(s)}ds},
\end{equation}
where~$\BRone$ is the ground state of~\eqref{eq:radial_CBNLS} with~$\sigma=4$
and~$d=1$, and that the blowup rate is slightly faster than a quartic root.
In other words, the blowup rate and profile are the same as those of peak-type
singular solutions of the one-dimensional critical BNLS,
see~\eqref{eq:intro_peak_SCBNLS}.

In Section~\ref{sec:shrinking}, we consider the case~$4/d<\sigma<4$, for
which~$0<\alphab<1$, see~\eqref{eq:psi_QN_alpha_BNLS}.
From power conservation we deduce that~$\alpha\ge\alphab$.
By analogy with the NLS, we expect that~$\alpha=\alphab$.
Therefore, the ring radius~$\rmax(t)\sim r_0 L^\alpha(t)$ decays to zero, but at a
slower rate than~$L(t)$.
%
%Informal asymptotic analysis shows that the blowup rate blowup rate 
%with~$p\le\frac1{3+\alpha}.$
By analogy with the NLS, we also expect that the blowup rate of these ring
solutions is given by~\eqref{eq:blowup_rate_with_p} with
\[
    p=\frac1{4-\frac{\sigma d-4}{\sigma(d-1)}}
    =\frac1{3+\alphab}.
\]
Numerical experiments support these predictions.

In Section~\ref{sec:critical}, we consider the critical BNLS ($\sigma=4/d$),
which corresponds to~$\alphab=1$.
Since the singular part of the solution has to be self-similar in~$r/L$,
see Theorem~\ref{thrm:self-similarity}, $\alpha$ must be equal to unity.
By the analogy with the NLS, the blowup rate is conjectured to be~$1/4$.
%the rigorous quartic-root lower bound, see Theorem~\ref{thrm:low-bound}.

\begin{figure}[h]
    \begin{center}
        \scalebox{0.7}{
        \input{BNLS_ring_class.pstex_t}}
        \mycaption{\label{fig:BNLS_class_diagram}
            Classification of singular ring-type solutions of the BNLS as a
            function of~$\sigma$ and~$d$.
            A:~subcritical case (no singularity).
            B: critical case, with equal-rate collapse.
            %, see Section~\ref{sec:critical}.
            C: $4/d<\sigma<4$, shrinking rings.
            %, see Section~\ref{sec:shrinking}.
            D:~$\sigma=4$, standing rings. %, see Section~\ref{sec:standing}.
            E:~$\sigma>4$, standing rings~\cite{Baruch_Fibich_Gavish:2009}.
        }
    \end{center}
\end{figure}
In summary, the BNLS singular ring-type solutions can be classified as follows
(see Figure~\ref{fig:BNLS_class_diagram}):
\begin{enumerate}   \renewcommand{\theenumi}{\Alph{enumi}}
    \item In the subcritical case~($\sigma d<4$), all BNLS solutions exist
        globally, hence no collapsing ring solutions exist.
    \item The critical case~$\sigma d=4$ corresponds to~$\alphab=1$ (equal-rate
        collapse).
        The blowup rate is~$p=1/4$.
    \item The supercritical case~$4/d<\sigma<4$ corresponds to~$0<\alphab<1$,
        hence the ring radius~$\rmax(t)$ decays to zero, but at a slower rate
        than~$L(t)$.
        The blowup rate is~$
            p %=\frac1{3+\alphab}
            =\frac1{4-\frac{\sigma d-4}{\sigma(d-1)}}
        $.
    \item The case $\sigma=4$ corresponds to~$\alphab=0$.
        Hence the solution is a singular standing ring.
        The self-similar profile~$\BQ$ is equal to that of peak-type solutions
        of the~$1D$ critical BNLS, and the blowup rate is slightly
        above~$p=1/4$.
    \item The case~$\sigma>4$ was studied
        in~\cite{Baruch_Fibich_Gavish:2009}.
        In this case, the solutions are of the standing-ring type, the
        self-similar profile is equal to that of peak-type solutions of
        the~$1D$ supercritical BNLS, and the blowup-rate is a quartic-root.
\end{enumerate}
\noindent
Thus, up to the change~$2\longrightarrow4$, this classification is, indeed,
completely analogous to that of singular ring-type solutions of the NLS 
(see Figure~\ref{fig:NLS_class_diagram}).

\subsection{Numerical Methodology}

The computations of singular BNLS solutions that focus by factors of~$10^8$
necessitated the usage of adaptive grids.
For our simulations we developed a modified version of the Static Grid
Redistribution method~\cite{Ren-00,SGR-08}, which is much easier to
implement in the biharmonic problem, and is easily extended to other evolution
equations, such as the nonlinear heat and biharmonic nonlinear heat
equations~\cite{Baruch_Fibich_Gavish:2009}.
The method of~\cite{SGR-08} also includes a mechanism for the prevention of
under-resolution in the non-singular region.
We extend this mechanism to prevent under-resolution in the transition layer
between the singular and non-singular regions.
See Section~\ref{sec:SGR} for further details.

\subsection{Critical exponents of singular ring solutions}\label{sec:discussion}

\begin{figure}[h]
    \begin{center}
    \includegraphics[clip,angle=-90,width=0.8\textwidth]%
        {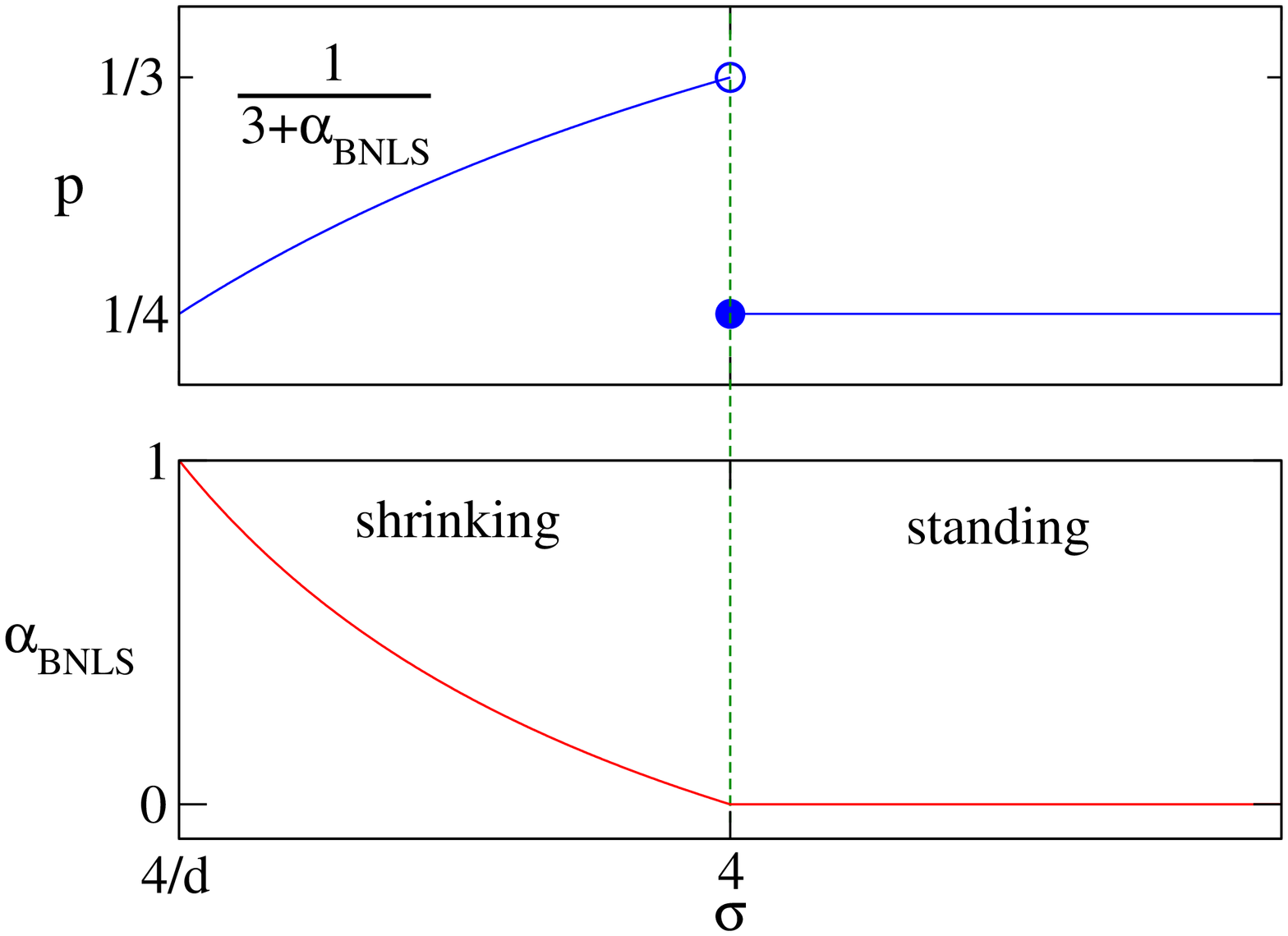}
    \mycaption{ \label{fig:ring_phase_transition}
        top: Blowup rate of singular ring solutions of the BNLS.
        The blowup rate increases monotonically from~$p=1/4$ at~$\sigma=4/d$
        to~$p=(1/3)-$ at~$\sigma=4-$.
        For~$\sigma = 4$ (full circle)~$p=1/4$ (with a loglog correction?) and
        for~$\sigma>4$,~$p\equiv1/4$.
        bottom: The shrinkage parameter~$\alphab$ of singular ring solutions of
        the BNLS.
        For~$4/d\le\sigma<4$,~$\alphab$ decreases monotonically from~$1$
        to~$0+$ (shrinking rings).
        For~$\sigma\ge4$,~$\alphab\equiv0$ (standing rings).
    }
    \end{center}
%D:\Research\My Work\Matlab\Super Critical Ring profile\IGR
\end{figure}
In Figure~\ref{fig:ring_phase_transition} (top) we plot the blowup rate~$p$ of
singular ring solutions of the BNLS, see~\eqref{eq:blowup_rate_with_p}.
As~$\sigma$ increases from~$4/d$ to~$4-$,~$p$ increases monotonically
from~$\frac14$ to~$\frac13-$.
At~$\sigma=4$, the blowup rate is slightly faster than a quartic root,
i.e.,~$p\approx\frac14$.
Finally,~$p=1/4$~for~$\sigma>4$.
Since
\[
    \lim_{\sigma\to4-}p=1/3,\qquad\lim_{\sigma\to4+}p=\frac14,
\]
the blowup rate has a discontinuity at~$\sigma=4$.
%Surprisingly, the blowup rate is not monotonically-increasing with~$\sigma$.
%For example, a ring solution of the NLS with~$\sigma=1.8$ blows up faster than
%a ring solution of the NLS with~$\sigma=2.2$.

The above results show that
{\em $\sigma=4$ is a critical exponent of singular ring solutions of the
BNLS}.
Intuitively, this is because the blowup dynamics changes from a
shrinking-ring~$(\sigma<4)$ to a standing-ring~$(\sigma\ge4)$, see
Figure~\ref{fig:ring_phase_transition} (bottom).
We can understand why~$\sigma=4$ is a critical exponent using the
following argument.
Standing-ring solutions are `equivalent' to singular peak solutions
of the one-dimensional NLS with the same nonlinearity exponent~$\sigma$
\cite{Baruch_Fibich_Gavish:2009}.
Since~$\sigma=4$ is the critical exponent for singularity formation in the
one-dimensional NLS, it is also the critical exponent for standing-ring blowup.
An analogous picture exists for the NLS, wherein the phase transition between
standing and shrinking rings occurs
at~$\sigma=2$~\cite{Baruch_Fibich_Gavish:2009}.

\section{\label{sec:standing}Singular standing rings ($\sigma=4$)}

In what follows, we show that collapse of ring-type singular solutions of
the BNLS with~$\sigma=4$ is ``the same'' as collapse of peak-type singular
solutions of the one-dimensional critical BNLS.

\subsection{\label{ssec:standing_analysis}Informal Analysis}

We consider ring-type singular solutions of the supercritical
BNLS~\eqref{eq:radial_BNLS} with~$\sigma=4$, that undergo a quasi self-similar
collapse with the asymptotic profile
\begin{equation}    \label{eq:shrinking-SS-s4}
        \psi_{\BQ}(t,r) =
            \frac{1}{L^{1/2}(t)}
            \BQ\left(\rho\right) e^{i \int^t\frac{1}{L^4(s)}ds}
            \qquad \rho=\frac{r-\rmax(t)}{L} .
            %\qquad -\rho_c<\rho<\rho_c.
\end{equation}
Here and throughout this paper, by quasi self-similar we mean
that~$\psi\sim\psi_{\BQ}$ in the singular ring
region~$\left.r-\rmax=\mathcal{O}(L)\right.$, or~$\rho=\mathcal{O}(1)$, but not
for~$0\le r<\infty$.

The asymptotic profile~\eqref{eq:shrinking-SS-s4} describes a standing ring
if~$
    \displaystyle \lim_{t\to\TCrit} \rmax(t) > 0
$.
We expect ring-type singular solutions of the BNLS with~$\sigma=4$ to collapse
as standing rings, for the following two reasons:
\begin{enumerate}
    \item By continuity, since ring-type singular solutions of the BNLS
        with~$\sigma>4$ are standing rings~\cite{Baruch_Fibich_Gavish:2009}.
    \item By analogy with singular ring-type solutions of the NLS with~$\sigma=2$,
        which are standing rings~\cite{SC_rings-07,Raphael-06,Raphael-08}.
\end{enumerate}

In the ring region~$r-\rmax=\mathcal{O}(L)$, as~$L\to0$, the terms of the
radial biharmonic operator~\eqref{eq:radial_bi_Laplacian} behave as \[
    \left[ \frac1{r^{4-k}}\partial_r^k\psi  \right]
    \sim \frac{[\psi]}{L^{k}}
    ,\qquad k=0,\dots,4.
\]
Therefore,~$\Delta_r^2\psi \sim \partial^4_r\psi$.
Hence, near the singularity, equation~\eqref{eq:radial_BNLS} reduces to
\[
    \psi(t,r) - \psi_{rrrr} + |\psi|^8\psi = 0,
\]
which is the one-dimensional critical BNLS.
Therefore, the singular solutions of the two equations are asymptotically
equivalent, i.e.,\[
    \psi_{\sigma=4,d}^{ring}(t,r)
    \sim
    \psi_{\sigma=4,d=1}^{peak}(t,x=r-\rmax(t)),
\]
where~$\psi_{\sigma=4,d=1}^{peak}$ is a peak-type solution of the
one-dimensional critical BNLS.

The above informal analysis suggests that the blowup dynamics of singular
standing-ring solutions of the BNLS~\eqref{eq:radial_BNLS} with~$d>1$
and~$\sigma=4$ is the same as the blowup dynamics of singular peak solutions of
the one-dimensional critical BNLS:
\begin{conj}    \label{conj:standing_ring_equals_peak}
    Let~$d>1$ and~$\sigma=4$, and let~$\psi$ be a singular ring-type solution
    of the BNLS~\eqref{eq:radial_BNLS}.
    Then,
    \begin{enumerate}
        \item The solution is a standing ring, i.e.,
            $ \displaystyle \lim_{t\to\TCrit} \rmax(t) > 0$.
        \item In the ring region, the solution approaches the~$\psi_{\BQ}$
            self-similar profile, see~\eqref{eq:shrinking-SS-s4}.
        %   \begin{equation}
        %       \psi(t,r) \sim \psi_{R,1D}(t,r),
        %       \qquad
        %       \rmax(t)-\rho_cL(t) < r < \rmax(t)+\rho_cL(t).
        %   \end{equation}
            %i.e.,~$\psi\sim \psi_{\BQ}~$ for~$r-\rmax=\mathcal{O}(L)$,
            %where~$\psi_{R,1D}$ is given by~\eqref{eq:shrinking-SS-s4}.
        \item The self-similar profile~$\psi_{\BQ}$ is given by
            \begin{equation}\label{eq:psiB}
                \psi_{\BQ}(t,r)= \psi_{R_{B,1D}}\left( t,x=r-\rmax(t) \right),
            \end{equation}
            where~$\psi_{\BRone}(t,x)$, see~\eqref{eqs:phi-QSS}, is the
            asymptotic profile of the one-dimensional critical BNLS.
        \item Specifically,~$\BQ(\rho)=\BRone(\xi)$,
            where~$\BRone(\xi)$ is the ground state
            of~\eqref{eq:standing-1D}.
        \item The blowup rate of~$L(t)$ is slightly faster than a quartic root,
            see~\eqref{eq:rate_slightly_mad}.
    \end{enumerate}
\end{conj}

\noindent In Section~\ref{ssec:standing_simulations} we provide numerical
evidence in support of Conjecture~\ref{conj:standing_ring_equals_peak}.

\subsection{\label{ssec:standing_simulations}Simulations}

\begin{figure}
    \centering
    \subfloat{%
        \label{fig:standing_ring_rmax}%
        \includegraphics[angle=-90,clip,width=\myfw\textwidth]%
            {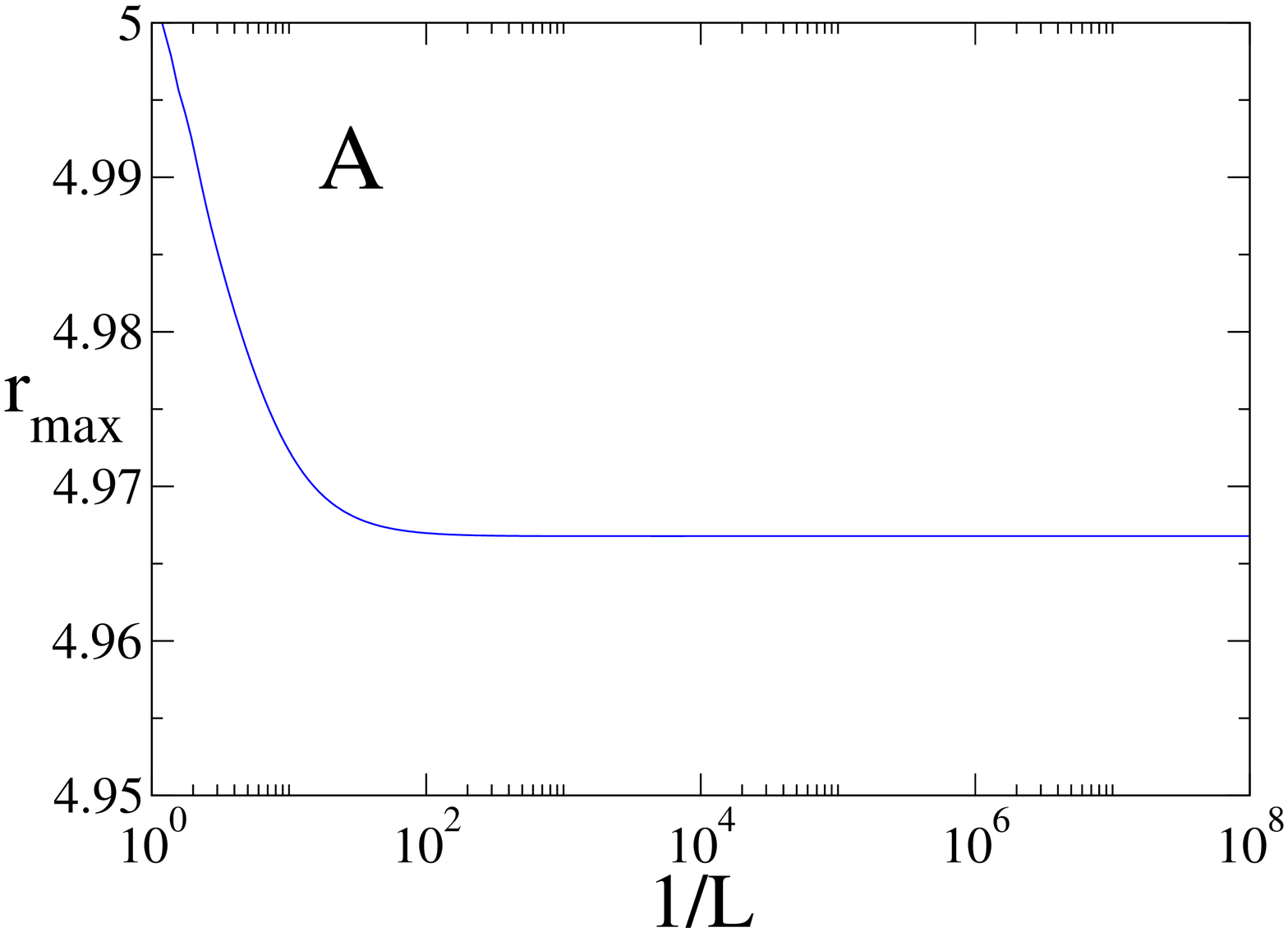}%
    }
    \subfloat{\label{fig:standing_ring_rescaled}%
        \includegraphics[angle=-90,clip,width=\myfw\textwidth]%
            {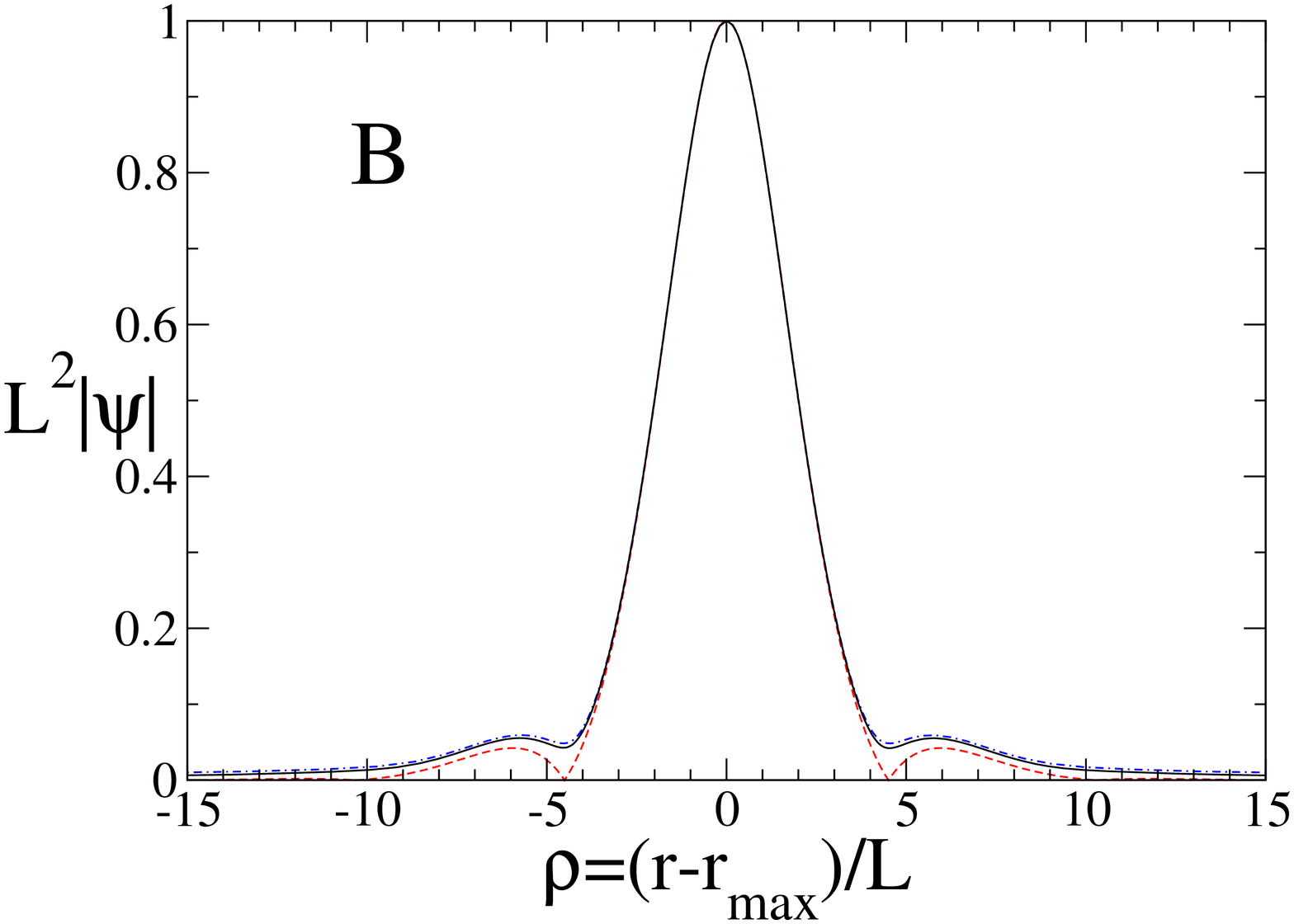}%
    }

    \mycaption{\label{fig:standing_ring_rmax_rescaled}
    A singular standing-ring solution of the supercritical
    BNLS~\eqref{eq:radial_BNLS} with~$d=2$ and $\sigma=4$.
    A) Ring radius~$\rmax$~as a function of the focusing level~$1/L$.
    B) The rescaled solution, see~\eqref{eq:psi_rescaled},
    at~$L(t)=10^{-4}$ (blue dash-dotted line) and~$L(t)=10^{-8}$ (black solid line).
    The two curves are indistinguishable.
    Red dashed line is the rescaled one-dimensional ground
    state~$\abs{\BRone(x)}$.
    }
\end{figure}
\begin{figure}
    \centering
    \subfloat{\label{fig:standing_ring_L}%
        \includegraphics[angle=-90,clip,width=\myfw\textwidth]%
            {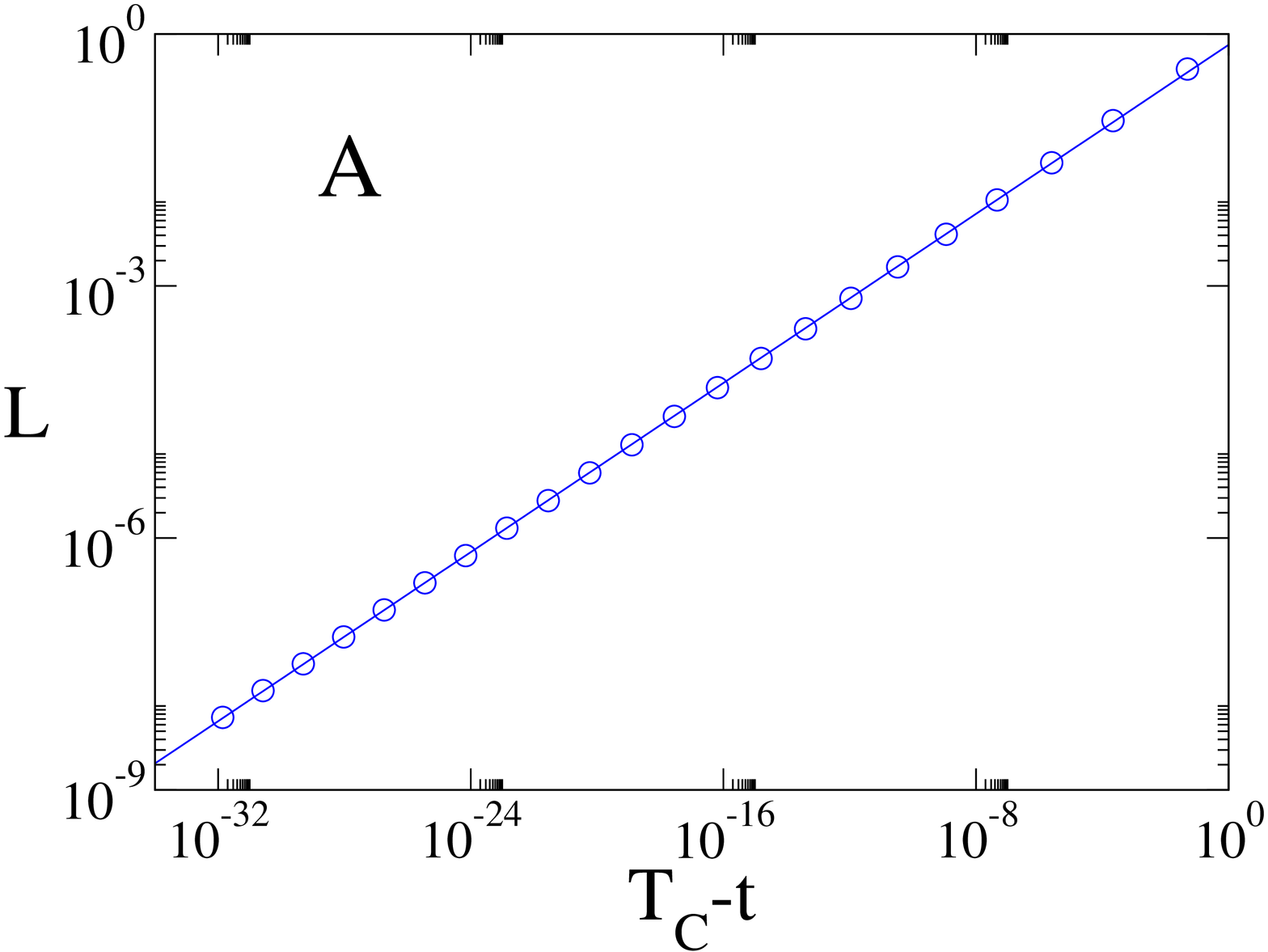}%
    }
    \subfloat{\label{fig:standing_ring_LtL3}%
        \includegraphics[angle=-90,clip,width=\myfw\textwidth]%
            {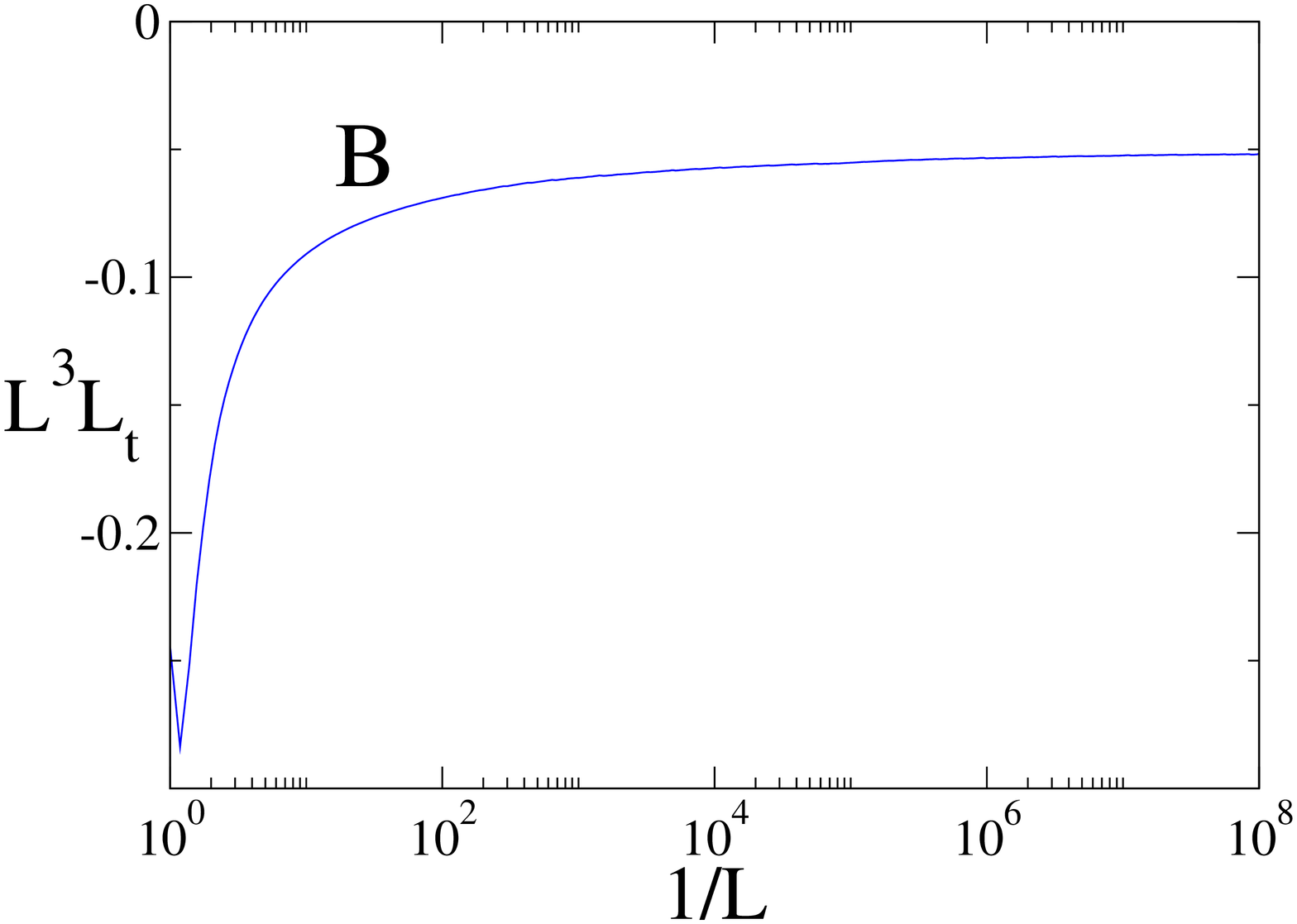}%
    }

    \mycaption{\label{fig:standing_ring_L_LtL3}
    Blowup rate of the solution of Figure~\ref{fig:standing_ring_rmax_rescaled}.
    A)~$L$ as a function of~$\left(\TCrit-t\right)$ on a
    logarithmic scale (circles).
    Solid line is~$\left.L=0.774(\TCrit-t)^{0.2523}\right.$.
    B)~$L_tL^3$ as a function of~$1/L$.
    }
\end{figure}

The radially-symmetric BNLS~\eqref{eq:radial_BNLS} with~$d=2$ and~$\sigma=4$
was solved with the initial
condition~$\left.\psi_0(r) = 2 e^{-(r-5)^2}\right.$.
The simulation was run up to~$L=\mathcal{O}(10^{-8})$.
Similar results were obtained with~$d=3$ and~$\sigma=4$ (data not shown).

We next test each item of
Conjecture~\ref{conj:standing_ring_equals_peak} numerically:
\begin{enumerate}
    \item The position of maximal amplitude~$
            \rmax(t)=\arg\max_r|\psi|
        $ approaches a positive constant as~$L\to 0$, see
        Figure~\ref{fig:standing_ring_rmax}, indicating that the solution
        collapses as a standing ring.
    \item The solution profiles, at the focusing levels of~$L=10^{-4}$
        and~$L=10^{-8}$, rescaled according to
        \begin{equation}    \label{eq:psi_rescaled}
            \psi_{\text{rescaled}}(t,\rho) =
                L^{2/\sigma}(t)\psi(t,\rmax(t)+\rho\cdot L),
            \qquad L(t)=\norm{\psi}_{\infty}^{-\sigma/2},
            %\qquad \rmax(t)=\arg\max_r|\psi|,
        \end{equation}
        %\]
        are almost indistinguishable, see Figure~\ref{fig:standing_ring_rescaled},
        indicating that the collapsing core is self-similar according
        to~\eqref{eq:shrinking-SS-s4}.
    \item Figure~\ref{fig:standing_ring_rescaled} also shows that the self-similar
        profile of the standing-ring solution is given by~$\BRone(\xi)$,
        the one-dimensional ground-state of equation~\eqref{eq:standing-1D}.
    \item To calculate the blowup rate of $\psi$, we first assume that~$
        %\begin{equation}\label{eq:L_assumption}
            L(t)\sim \kappa(\TCrit-t)^p,
        %\end{equation}
        $ and find the best fitting~$\kappa$ and~$p$, see
        Figure~\ref{fig:standing_ring_L}.
        In this case~$p \approx 0.2523$, indicating that the blowup rate is a
        quartic-root or slightly faster.
    \item In order to check whether the blowup rate of~$L$ is slightly faster
        than a quartic-root, we compute the
        limit~$\displaystyle{\lim_{t\to\TCrit}}L^3L_t$.
        Recall that for a quartic-root blowup rate~$
            L(t)\sim\kappa(\TCrit-t)^{1/4}
        $ with~$\kappa>0$,
        \[
            \lim_{t\to \TCrit}L^3L_t=-\frac{\kappa^4}4<0,
        \]
        while for a faster-than-a quartic-root blowup rate,
        see~\eqref{eq:rate_slightly_mad}, $L^3L_t$ goes to zero.
        Figure~\ref{fig:standing_ring_LtL3} shows that $L^3L_t$ does not
        approach a negative constant, but increases slowly towards~$0^-$, implying
        that the blowup rate is slightly faster than a quartic root.
\end{enumerate}
\noindent
Note that the initial condition~$\psi_0=2e^{-(r-5)^2}$ is quite different from the
asymptotic profile~$\psi_{\BQ}$, indicating that the standing-ring~$\psi_{\BQ}$
profile~\eqref{eq:psiB} is an attractor in the radial case.

\section{\label{sec:shrinking}Shrinking-ring solutions of the supercritical BNLS ($4/d<\sigma<4$)}
%\section{\label{sec:ring} Ring-type singular solutions}

In this section, we consider the regime~$4/d<\sigma<4$.
In the NLS analogue ($2/d<\sigma<2$), the asymptotic profile has a
``double-lens'' radial-phase term, see~\eqref{eq:psi_QN_profile}, whose explicit
form is used in the asymptotic calculation of the blowup rate and shrinking
rate.
In contrast, for the BNLS we do not know the corresponding ``double-lens''
radial phase term, but only the amplitude~$\abs{\psi_{\BQ}}$.
Therefore, the results of the asymptotic analysis are weaker, and we need to 
rely on the analogy between the NLS and the BNLS.

\subsection{\label{ssec:shrinking_analysis}Informal analysis}

We consider singular ring-type solutions of the supercritical BNLS
with~$4/d<\sigma<4$, that undergo a quasi self-similar collapse with the
asymptotic profile~$\psi_{\BQ}$, whose amplitude is given by
\begin{equation}    \label{eq:shrinking-SS}
    %\begin{gathered}
        \abs{\psi_{\BQ}(t,r)} =
            \frac{1}{L^{2/\sigma}(t)}
            \abs{\BQ(\rho)},
        %e^{i\tau(t)}, \qquad \tau(t) = \int^{t}\frac{1}{L^4(s)}ds \\
        \quad \rho=\frac{r-\rmax(t)}{L},
        %\qquad -\rho_c\le \rho \le \rho_c,
        \quad \rmax(t) \sim r_0 L^\alpha(t).
        %\qquad \alpha\le1.
    %\end{gathered}
\end{equation}
As before, we assume that~$\psi\sim\psi_{\BQ}$ in the
region~$\rmax-r=\mathcal{O}(L)$, i.e., for~$|\rho|\le\rho_c=\mathcal{O}(1)$.
We assume that~$\alpha\le1$, since otherwise the rings are unstable.
Indeed, if $\alpha>1$, then~$\rho = r/L + o(1)$, and the rings
eventually evolve into a peak solution.

We first derive a lower bound for~$\alpha$:
\begin{lem} \label{lem:shrinking_alpha_bound}
    Let~$4/d<\sigma<4$, and let~$\psi$ be a ring-type singular solution of the
    BNLS equation~\eqref{eq:radial_BNLS}, whose asymptotic profile is of the
    form~\eqref{eq:shrinking-SS} with~$\alpha\le1$.
    Then,
    \[
        \alpha \ge \alphab,
    \]
    where
    \[
        \alphab = \frac{4-\sigma}{\sigma(d-1)}>0\,.
    \]
    Therefore, the ring is shrinking, i.e.,~$
        \displaystyle \lim_{t\to\TCrit}\rmax(t)=0
    $.
    \begin{proof}
        First, since~$4/d<\sigma<4$, then~$0<\alphab<1$.
        %Therefore,~$\rmax = r_0L^\alpha\to0$ as~$L\to 0$.
        The power of the collapsing core~$\psi_{\BQ}$ is
        \begin{eqnarray*}
            \norm{\psi_{\BQ}}_2^2
            &=&
            L^{-4/\sigma}\int_{r=\rmax-\rho_c\cdot L(t)}^{\rmax+\rho_c\cdot L(t)}
            \abs{\BQ\left( \frac{r-\rmax}{L} \right)}^2r^{d-1}dr \\
            &=&
            L^{-4/\sigma}\int_{\rho=-\rho_c}^{\rho_c}
            \abs{\BQ(\rho)}^2(L\rho+r_0 L^\alpha)^{d-1}(Ld\rho).
        \end{eqnarray*}
        In the case~$\alpha<1$, we have that~$
            L|\rho|\le L\rho_c \ll r_0L^\alpha
        $, hence~$L\rho+r_0 L^\alpha\sim r_0 L^\alpha$.
        Therefore, \[
        \norm{\psi_{\BQ}}_2^2 \sim L^{1-4/\sigma+\alpha(d-1)}(t)
            \cdot r_0^{\alpha(d-1)}
            \int_{\rho=-\rho_c}^{\rho_c} \abs{\BQ(\rho)}^2d\rho.
        \]
        In the case~$\alpha=1$, we have that~$
            L\rho+r_0 L^\alpha = (r_0+\rho) L
        $, hence \[
        \norm{\psi_{\BQ}}_2^2 \sim
            L^{1-4/\sigma+\alpha(d-1)}(t)
            \int_{\rho=-\rho_c}^{\rho_c} \abs{\BQ(\rho)}^2
            \left( r_0+\rho \right)^{d-1} d\rho.
        \]
        In both cases~$
            \norm{\psi_{\BQ}}_2^2 = \mathcal{O}\left(
                L^{1-4/\sigma+\alpha(d-1)}
            \right)
        $.
        Since~$
            \norm{\psi_{\BQ}}_2^2
                \le \norm{\psi}_2^2
                = \norm{\psi_0}_2^2<\infty
        $, then~$L^{1-4/\sigma+\alpha(d-1)}$ has to be bounded as~$L\to0$.
        Therefore,~$
            1-4/\sigma+\alpha(d-1)  \ge 0
        $, from which the result follows.
    \end{proof}
\end{lem}

Let
\[
    P_{\text{collapse}}=
        \liminf_{\varepsilon\to 0+}
        \lim_{t\to \TCrit}
        \int_{r<\varepsilon}|\psi|^2r^{d-1}dr
\]
be the amount of power that collapses into the singularity.
We say that the solution~$\psi$ undergoes a {\em strong collapse}
if~$P_{\text{collapse}}>0$, and a {\em weak collapse} if~$P_{\rm collapse}=0$.

\begin{cor} \label{cor:shrinking_ring_strength}
    Under the conditions of Lemma~\ref{lem:shrinking_alpha_bound},~$\psi_{\BQ}$
    undergoes a strong collapse if~$\alpha=\alphab$, and a weak collapse
    if~$\alpha>\alphab$.
\end{cor}
\begin{proof}
    This follows directly from the proof of Lemma~\ref{lem:shrinking_alpha_bound}.
\end{proof}

In the NLS with~$2/d<\sigma<2$, the shrinking rings undergo a strong collapse
with~$\alpha=\alphan$, see~\eqref{eq:psi_QN_alphan}.
Therefore, by analogy, we expect that the shrinking rings of the BNLS will also
undergo a strong collapse, in which case~$\alpha=\alphab$.

The blowup rate of singular shrinking-ring solutions of the
NLS with~$2/d<\sigma<2$ is~\cite{SC_rings-07}
\[
    p = \frac{1}{1+\alphan}
    = \frac1{2-\frac{\sigma d-2}{\sigma(d-1)}}.
\]
From the analogy of the BNLS with the NLS (up to the change $2\to4$), we
expect that the blowup rate of singular shrinking rings of the BNLS is
\[
    p = \frac1{4-\frac{\sigma d-4}{\sigma(d-1)}} = \frac{1}{3+\alphab}.
\]
Therefore, we have the following Conjecture:
\begin{conj}    \label{conj:shrinking_ring}
    Let~$d>1$ and~$4/d<\sigma<4$, and let~$\psi$ be a singular ring-type
    solution of the BNLS~\eqref{eq:radial_BNLS}.
    Then,
    \begin{enumerate}
        \item The solution is quasi self-similar, i.e.,~$\psi\sim\psi_{\BQ}$
            for~$r-\rmax = \mathcal{O}(L)$, where~$\psi_{\BQ}$ is given
            by~\eqref{eq:shrinking-SS}.
        \item The solution is a shrinking ring, i.e.,~$
            \displaystyle \lim_{t\to\TCrit} \rmax(t) = 0
        $.
        \item The shrinking rate is
            \begin{equation}    \label{eq:alpha_BNLS}
                \alpha = \alphab = \frac{4-\sigma}{\sigma(d-1)}.
            \end{equation}
        Specifically,~$0<\alpha<1$.
        \item The blowup rate is
            \begin{equation}    \label{eq:shrinking_rate}
                L(t)\sim\kappa( \TCrit-t )^p,
                \qquad\quad
                p = \frac{1}{3+\alphab}
                =\frac{1}{
                    4-\frac{\sigma d-4}{\sigma(d-1)}
                }.
                %=\frac{
                %   \sigma(d-1)
                %}{
                %   3\sigma d +4(1-\sigma)
                %}.
            \end{equation}
            Specifically,~$\frac14<p<\frac13$.
%       \item The self-similar profile~$\BQ(\rho)$ is given by
%       \begin{equation}    \label{eq:shrinking_Bs}
%           -B
%           -i \frac{\kappa^{3+\alpha} (1-\alpha)r_0}{3+\alpha} B_\rho
%               - B_{\rho\rho\rho\rho} + \abs{B}^{2\sigma}B = 0,
%       \end{equation}
%       where~$\kappa>0$ is the coefficient in~\eqref{eq:shrinking_rate}.
%       Therefore,~$B$ is different from the ground-state solution of \[
%           -R(\xi)-R_{\xi\xi\xi\xi} + \abs{R}^{2\sigma}R = 0.
%       \]
    \end{enumerate}
\end{conj}

\noindent In Section~\ref{ssec:shrinking_simulations} we provide numerical
evidence in support of Conjecture~\ref{conj:shrinking_ring}.

\subsection{\label{ssec:shrinking_simulations}Simulations}

\begin{figure}
    \centering
    \subfloat{%
        \label{fig:shrinking_ring_rmax}%
        \includegraphics[angle=-90,clip,width=\myfw\textwidth]%
            {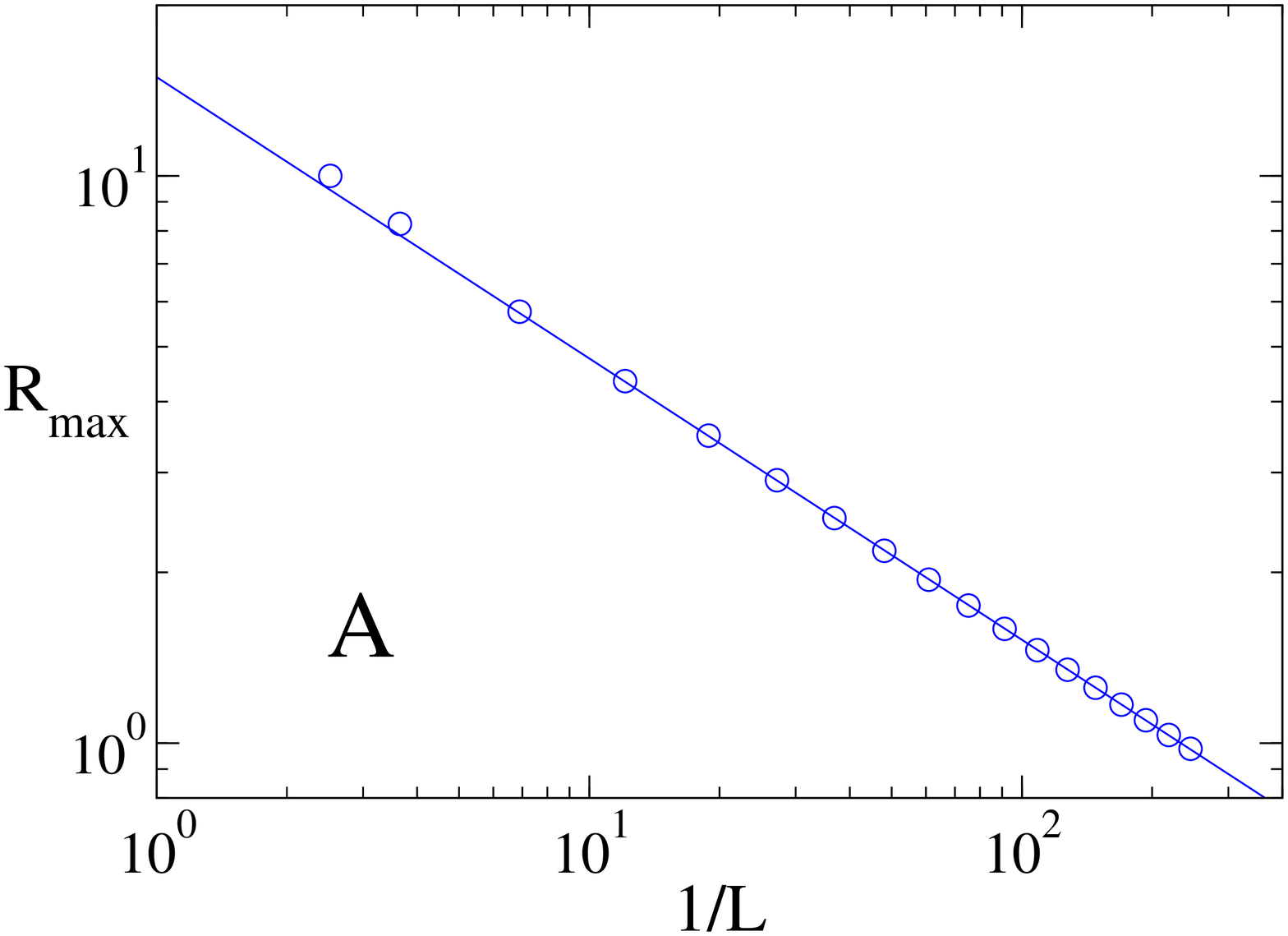}%
    }
    \subfloat{\label{fig:shrinking_ring_rescaled}%
        \includegraphics[angle=-90,clip,width=\myfw\textwidth]%
            {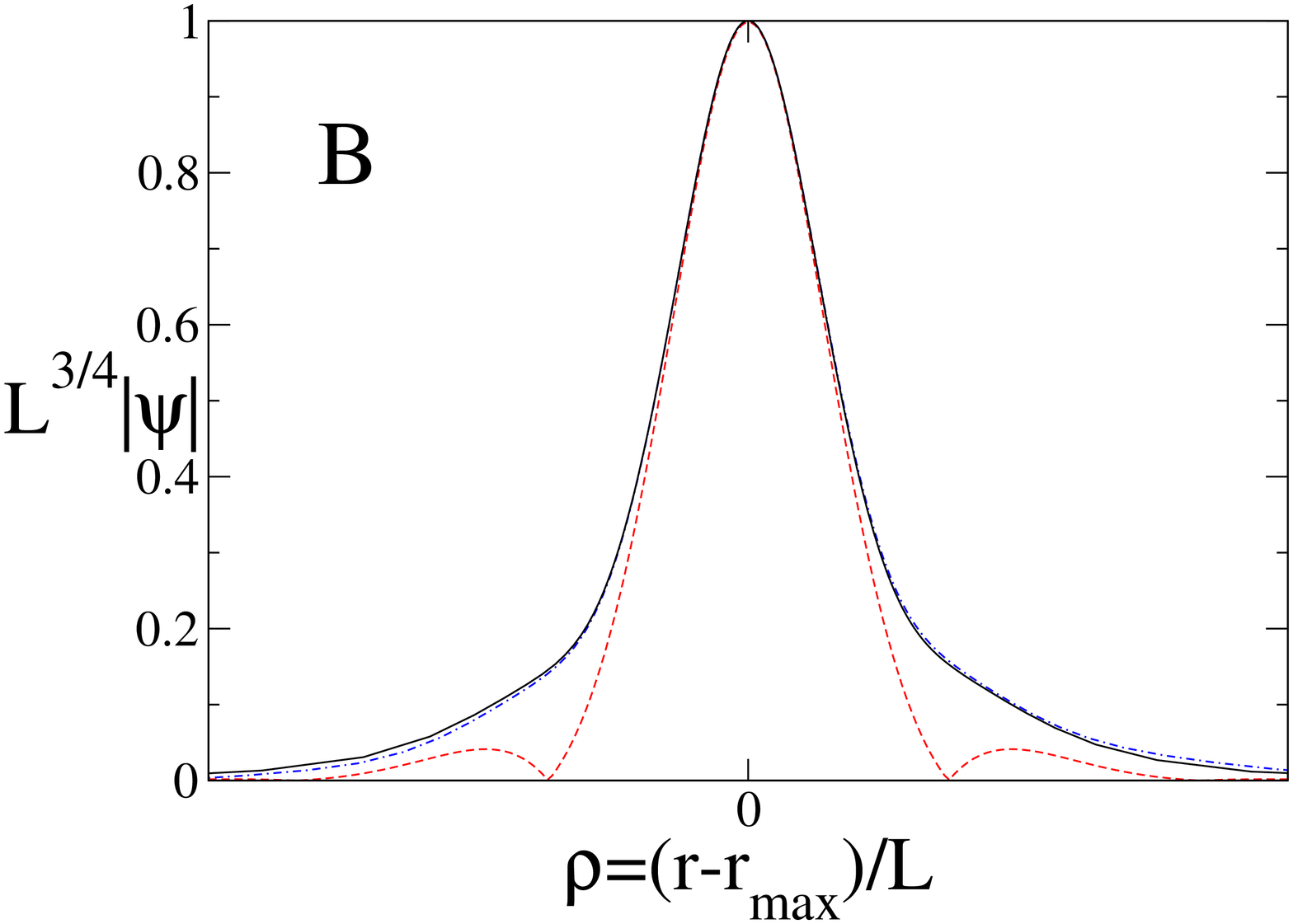}%
    }

    \mycaption{\label{fig:shrinking_ring_rmax_rescaled}
    Ring-type singular solution of the supercritical BNLS~\eqref{eq:radial_BNLS}
    with~$d=2$ and $\sigma=8/3$.
    A)~$\rmax$ as a function of the focusing factor~$1/L$.
    Solid line is~$\left.\rmax= 14.9 L^{0.496}\right.$.
    B) The rescaled solution, see~\eqref{eq:psi_rescaled},
    at~$L(t)=10^{-1}$ (blue dash-dotted line) and~$L(t)=10^{-2}$ (black solid line).
    %The curves are indistinguishable.
    Red dashed line is the rescaled one-dimensional ground
    state~$\abs{\BRone(x)}$.
    }
\end{figure}
\begin{figure}
    \centering
    \subfloat{\label{fig:shrinking_ring_L}%
        \includegraphics[angle=-90,clip,width=\myfw\textwidth]%
            {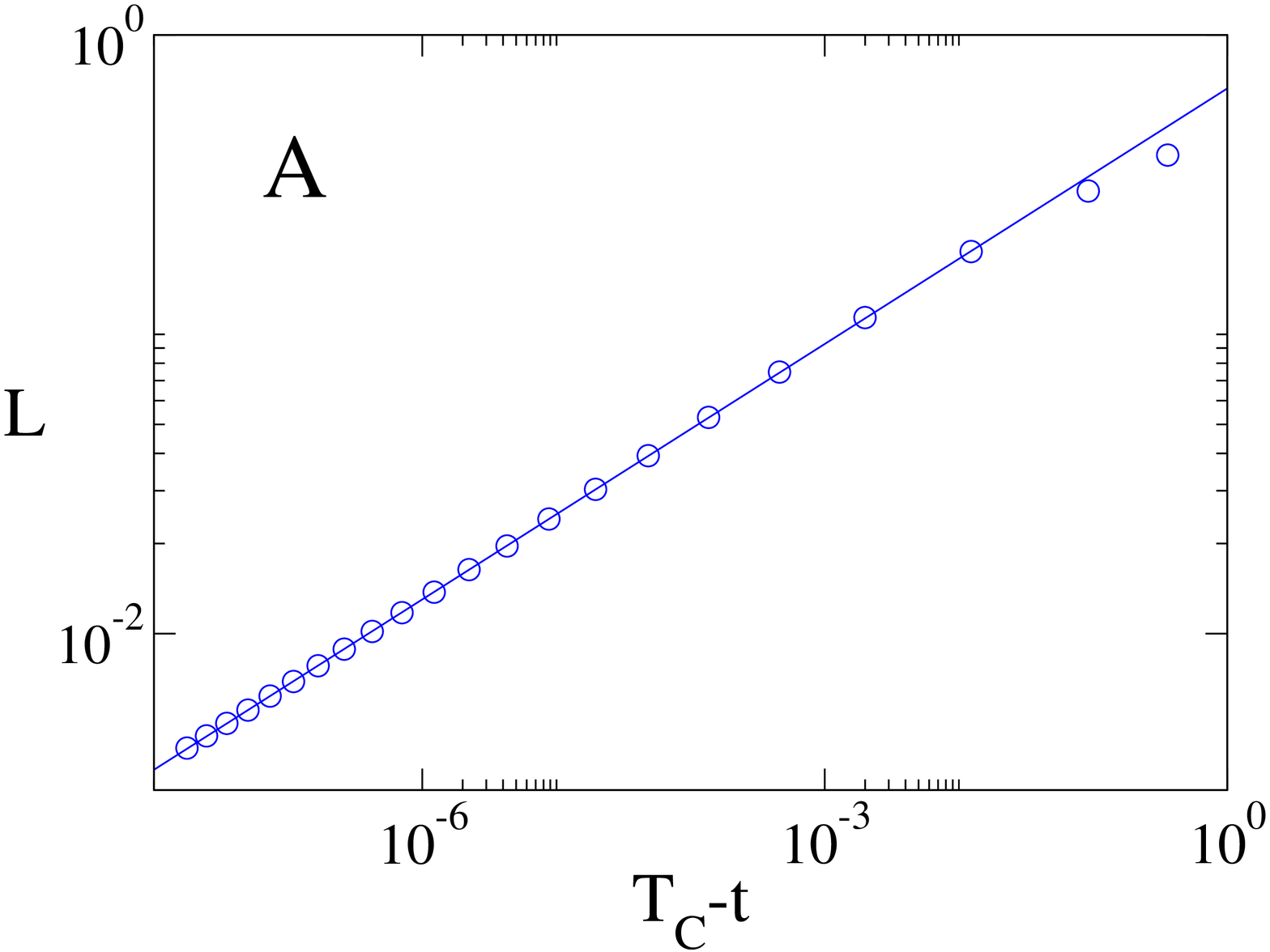}%
    }
    \subfloat{\label{fig:shrinking_ring_LtL2a}%
        \includegraphics[angle=-90,clip,width=\myfw\textwidth]%
            {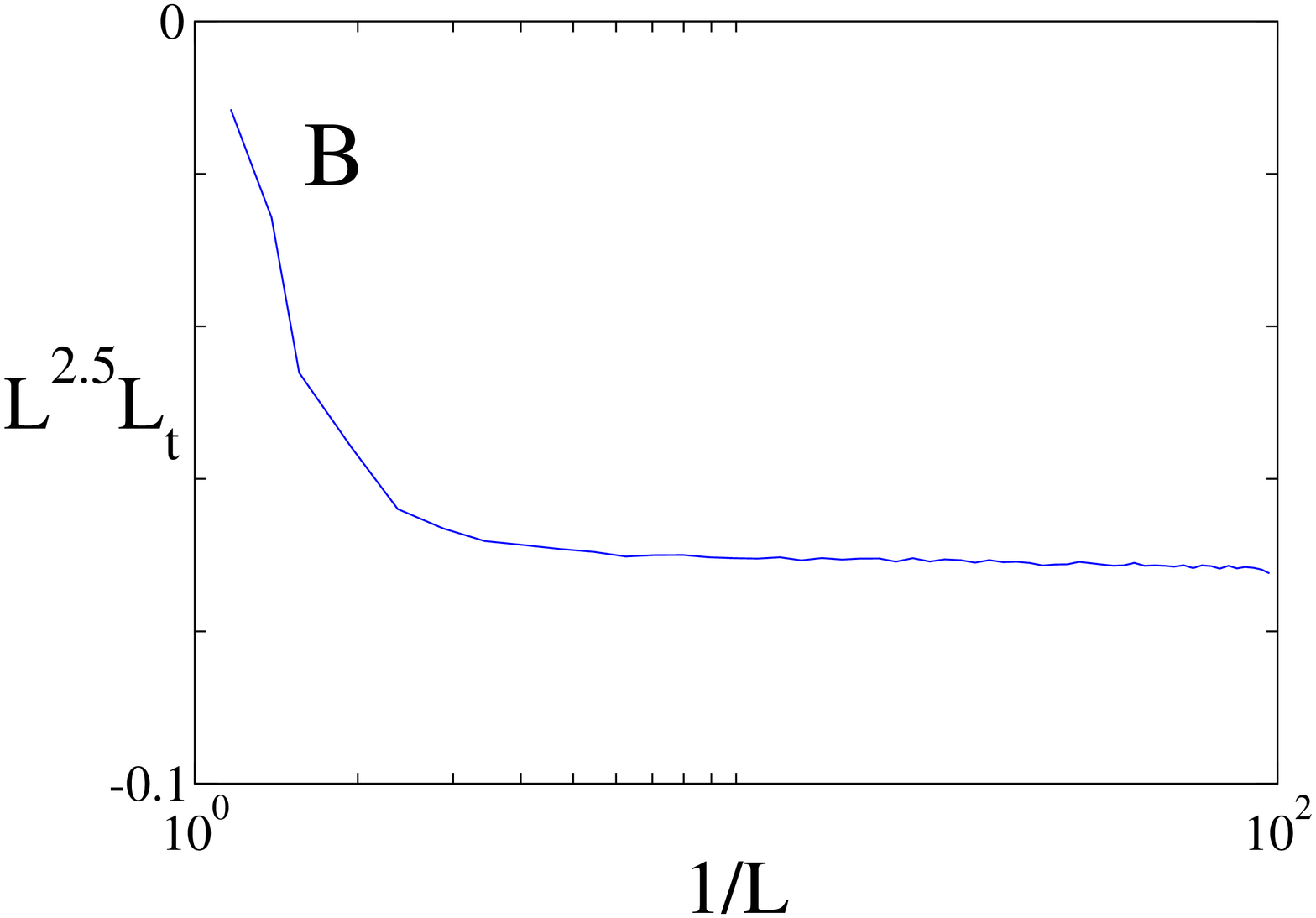}%
    }

    \mycaption{\label{fig:shrinking_ring_L_LtL3}
    Blowup rate of the solution of Figure~\ref{fig:shrinking_ring_rmax_rescaled}.
    A) Solution width~$L$ as a function of~$\left(\TCrit-t\right)$ on a
    logarithmic scale (circles).
    Solid line is~$\left.L=0.662(\TCrit-t)^{0.2844}\right.$.
    B)~$L^{2.5}L_t$ as a function of~$1/L$.
    }
\end{figure}

The supercritical BNLS equation with $d=2$ and $\sigma=8/3$ was solved with the
initial condition~$\psi_0= 2 e^{(r-10)^2}$.
The simulation was run up to a focusing level of~$L(t)=10^4$.

We next test Conjecture~\ref{conj:shrinking_ring} numerically, clause by
clause.
\begin{enumerate}
    \item Figure~\ref{fig:shrinking_ring_rmax} shows that the ring shrinks at a rate
        of~$\rmax(t)\sim14.9L^\alpha(t)$ with $\alpha\approx0.496$, which is close to
        the predicted value of \[
        \alphab
        %\left( \sigma={\textstyle\frac83},d=2 \right)
            = \frac{4-\frac83}{\frac83(2-1)}
            = \frac12\,.
        \]
        % c=?
    \item In Figure~\ref{fig:shrinking_ring_rescaled} we plot the solution, rescaled
        according to~\eqref{eq:psi_rescaled}, at the focusing levels~$1/L=10$
        and~$1/L=100$.
        The two curves are indistinguishable for~$\rho=\mathcal{O}(1)$, but not for
        all~$\rho$, showing that the solution undergoes a quasi self-similar
        collapse with the~$\psi_{\BQ}$ profile~\eqref{eq:shrinking-SS}.
    \item Figure~\ref{fig:shrinking_ring_rescaled} shows that the self-similar
        profile of the standing-ring solution is close to~$\BRone(\xi)$,
        the one-dimensional ground-state of equation~\eqref{eq:standing-1D}.
    \item Figure~\ref{fig:shrinking_ring_L} shows that~$
            L(t)\sim 0.662(\TCrit-t)^{0.2844}.
        $
        Therefore, the calculated blowup rate~$
            p=0.282
        $ is close to the predicted value of~$
            p = \frac{1}{3+\alphab} = 1/3.5
            \approx 0.2857
        $\,.
    \item In order to check whether the blowup rate of~$L$ is
        exactly~$p=1/3.5$, we compute the limit~$
            \displaystyle{\lim_{t\to \TCrit}}L^{2.5}L_t
        $.
        Recall that if~$
            L(t)\sim\kappa(\TCrit-t)^{1/3.5}
        $, then
        \[
            \lim_{t\to \TCrit}L^{2.5}L_t=-\frac{\kappa^{3.5}}{3.5}<0,
        \]
        while for a faster blowup rate,~$L^{2.5}L_t\to0$,
        and for a slower blowup rate,~$L^{2.5}L_t\to-\infty$.
        Figure~\ref{fig:shrinking_ring_LtL2a} shows that $L^{2.5}L_t$ converges
        to a negative constant, implying that the blowup rate is
        exactly~$p=1/3.5$.
\end{enumerate}

\begin{figure}
    \centering
    \subfloat{%
        \label{fig:shrinking_rings_rates_alpha_d2}%
        \includegraphics[angle=-90,clip,width=\myfwb\textwidth]%
            {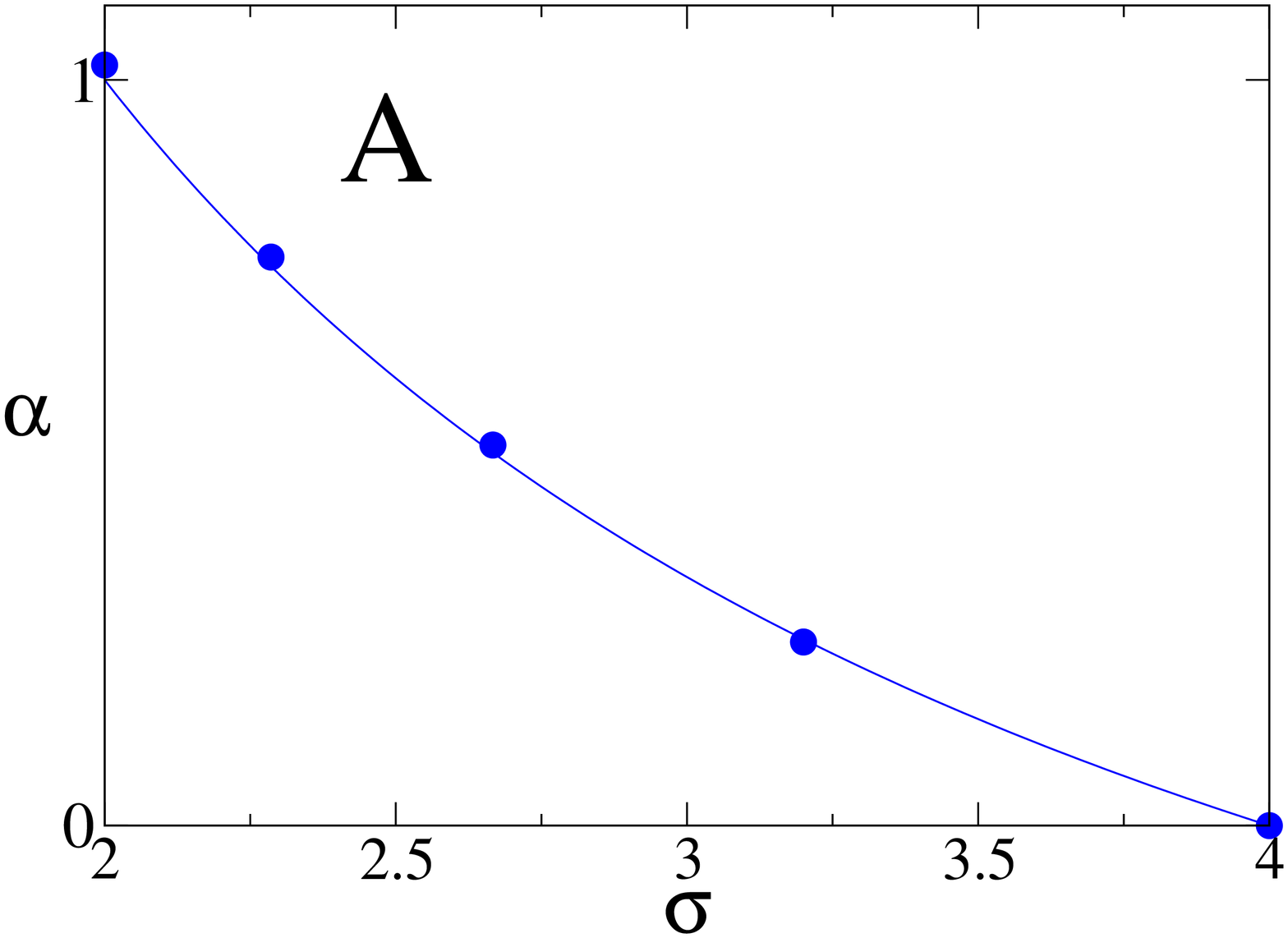}%
    }
    \subfloat{%
        \label{fig:shrinking_rings_rates_alpha_d3}%
        \includegraphics[angle=-90,clip,width=\myfwb\textwidth]%
            {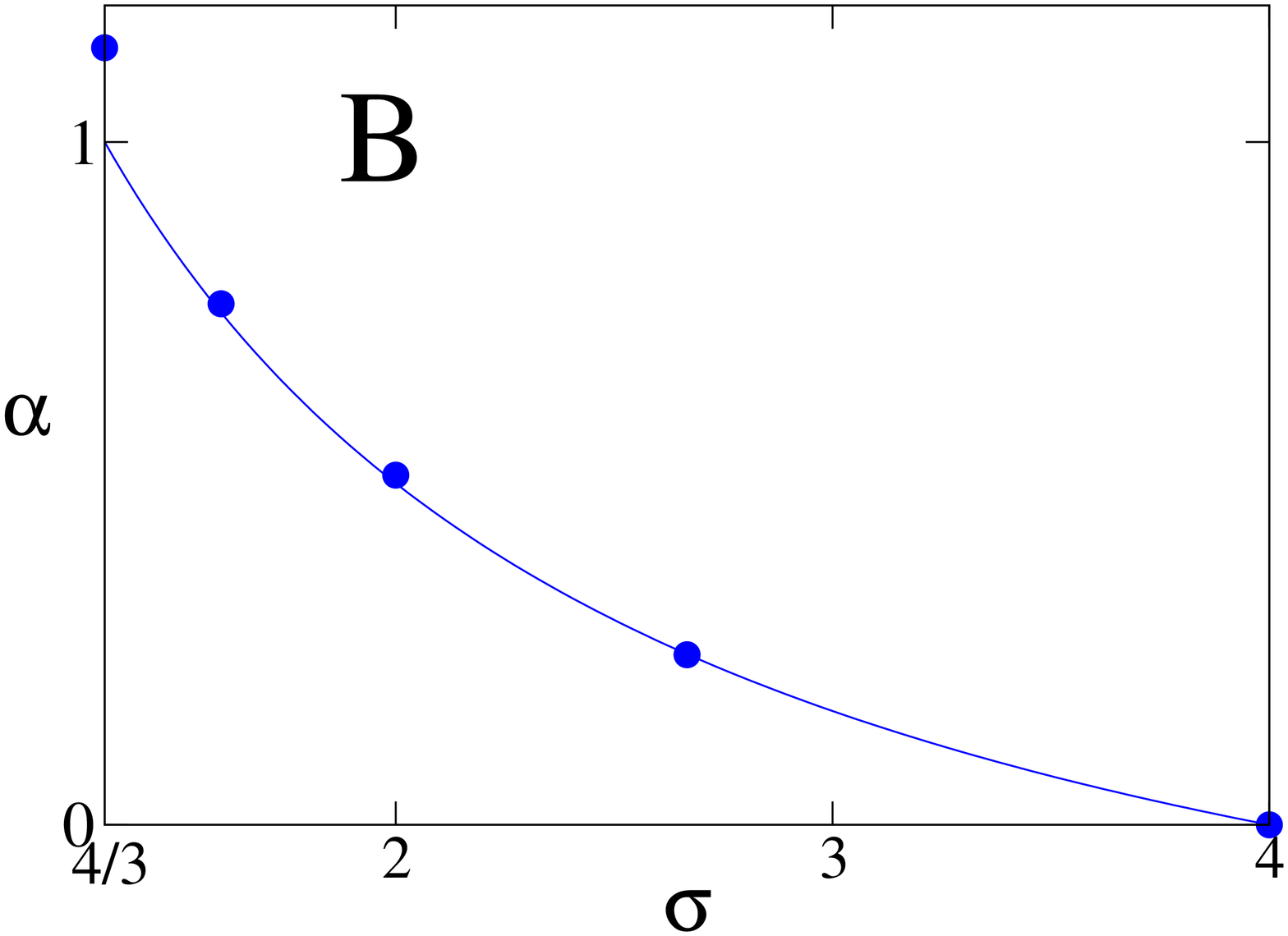}%
    }

    \mycaption{\label{fig:shrinking_rings_alpha}
    The numerical shrinking rate~$\alpha$
    (circles) for ring-type singular solutions of the BNLS
    with~$4/d\le\sigma\le4$.
    The solid line is~$\alpha=\alphab(\sigma,d)$,
    see~\eqref{eq:alpha_BNLS}.
    %The values of~$\sigma$ and~$d$ are given in Table~\ref{tab:shrinking_a_ps}.
    A)~$d=2$ and~$\sigma=2,16/7,8/3,16/5$~and $4$.
    B)~$d=3$ and~$\sigma=4/3,8/5,2,8/3$~and $4$.
    }
\end{figure}

\begin{figure}
    \centering
    \subfloat{%
        \label{fig:shrinking_rings_p_d2}%
        \includegraphics[angle=-90,clip,width=\myfwb\textwidth]%
            {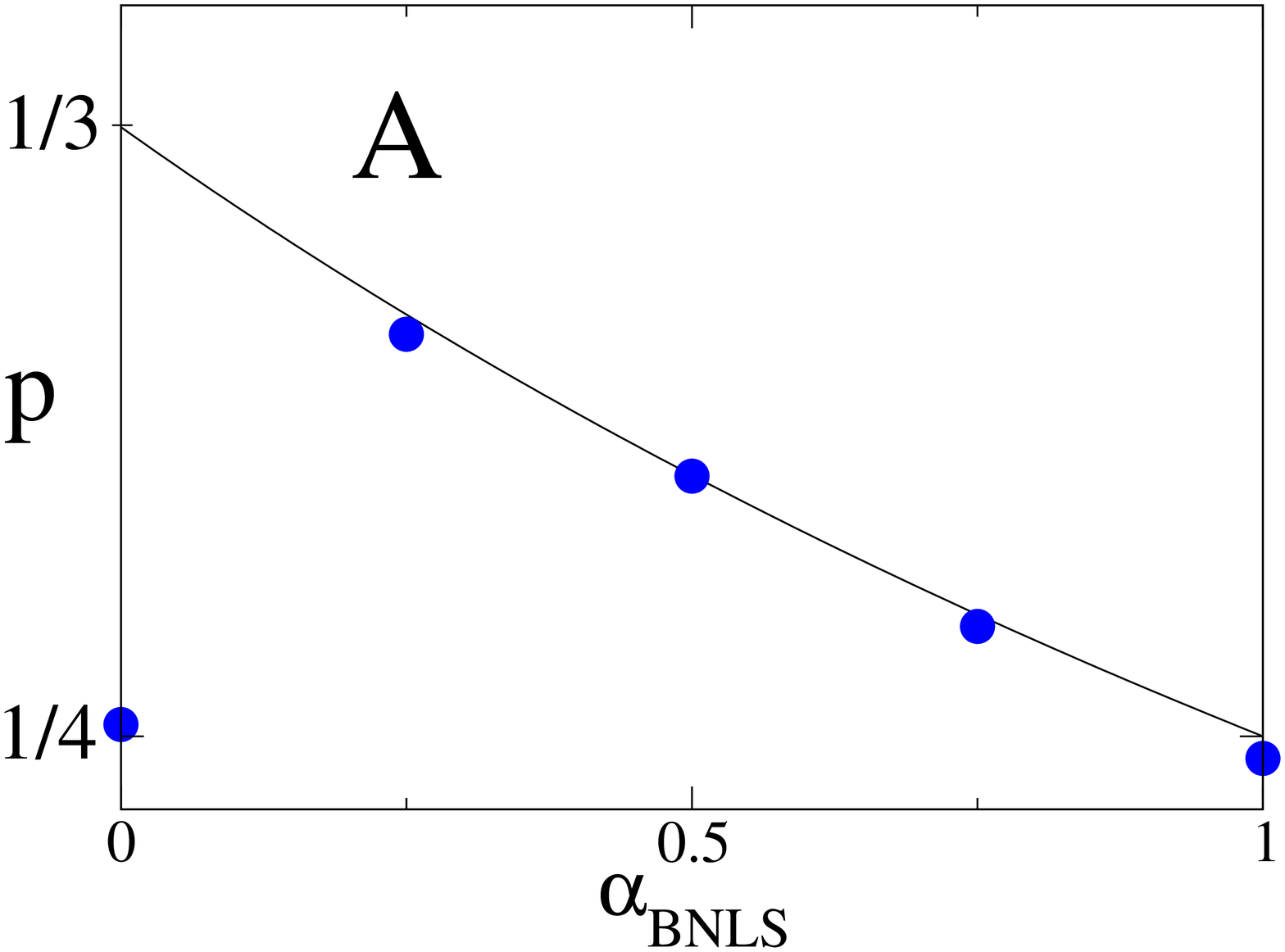}%
    }
    \subfloat{%
        \label{fig:shrinking_rings_p_d3}%
        \includegraphics[angle=-90,clip,width=\myfwb\textwidth]%
            {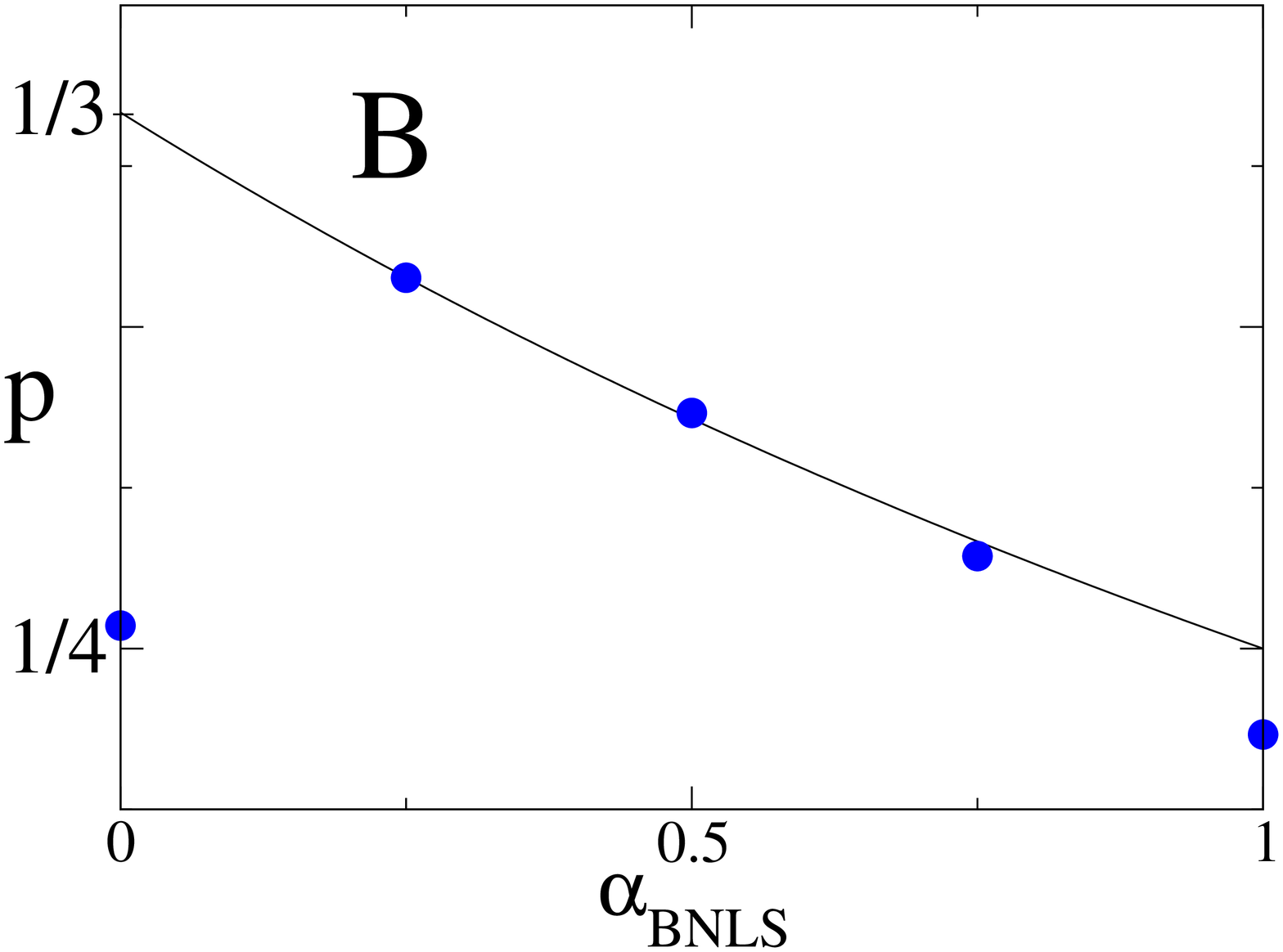}%
    }

    \mycaption{\label{fig:shrinking_rings_p}
        Same as Figure~\ref{fig:shrinking_rings_alpha}, for the numerical
        blowup rate~$p$, defined by~$L(t)\sim c(\TCrit-t)^p$, as a function
        of~$\alphab$.
        Solid line is~$p=\frac1{3+\alphab}$.
        The calculated values of~$p$ for~$1/2\le\alphab\le1$ are slightly lower
        than the predicted value~$\frac{1}{3+\alphab}$.
    }
\end{figure}

In Figure~\ref{fig:shrinking_rings_alpha} we present the numerical values
of the shrinking parameter~$\alpha$, defined by~$\rmax(t)\sim cL^\alpha$
for ten different values of~$(\sigma,d)$.
In all cases, the value of~$\alpha$ is very close to~$\alphab$,
see~\eqref{eq:alpha_BNLS}.
In Figure~\ref{fig:shrinking_rings_p} we present the numerical values of the
blowup rate~$p$ for the same simulations, and find that they are close
to~$p=1/\left( 3+\alphab \right)$ for~$\alphab>0$.
At~$\alphab=0+$, $p$ has a jump discontinuity to~$1/4$, in accordance
with Conjecture~\ref{conj:standing_ring_equals_peak}.
The discontinuity at~$\alphab=0$ ($\sigma=4$) is a manifestation of the
phase transition from shrinking to standing rings, see
Section~\ref{ssec:intro_summary}.

\section{\label{sec:critical}Equal-rate shrinking rings (critical BNLS)}

\subsection{\label{ssec:critical_analysis}Informal Analysis}

We consider singular ring-type solutions of the critical BNLS, that undergo a
quasi self-similar collapse with the asymptotic profile
\begin{equation}    \label{eq:crit-ring_SS0}
%\begin{gathered}
    \psi_{\BQ}(t,r) =
        \frac{1}{L^{d/2}(t)}
        \BQ \left(
            \frac{r-\rmax(t)}{L(t)}
        \right)
        e^{i \int^{t}\frac{1}{L^4(s)}ds}, %\\
    \quad \rmax \sim r_0 L^\alpha.
%\end{gathered}
\end{equation}
\begin{lem} \label{lem:crit-ring-alpha-1}
    Let~$\psi_{\BQ}(t,r)$, see~\eqref{eq:crit-ring_SS0}, be the asymptotic
    profile of singular ring-type solutions of the critical
    BNLS~\eqref{eq:radial_CBNLS}.
    Then,~$\alpha=1$.
    \begin{proof}
        Theorem~\ref{thrm:self-similarity} implies that the collapsing core of
        singular solutions of the critical BNLS is self-similar in~$r/L$, i.e., \[
            |\psi(t,r)|\sim \frac1{L^{d/2}}
                \left|\Psi\left( \frac rL \right)\right|.
        \]
        Therefore, a singular solution of the critical BNLS is ring-type if and
        only if~$\abs{\Psi(\rho)}$ attains its maximum at some~$\rho_{\max}>0$.
        Hence,~$\rmax(t) = \rho_{\max}\cdot L(t)$. Therefore,~$\alpha=1$.
    \end{proof}
\end{lem}

By Theorem~\ref{thrm:low-bound}, the blowup rate is lower-bounded by~$1/4$.
We recall that ring-type singular solutions of the critical NLS have a square-root
blowup rate~\cite{Gprofile-05}.
Therefore, we expect that ring-type singular solutions of the critical BNLS have a
quartic-root blowup rate.

In summary, we conjecture the following:
\begin{conj}    \label{conj:crit_ring_rate_profile}
        Let $\psi$ be ring-type singular solution of the critical
        BNLS~\eqref{eq:radial_CBNLS}.
    Then,
    \begin{subequations} \label{eqs:crit_ring_QSS}
        \begin{enumerate}
            \item The solution undergoes an equal-rate collapse,
                i.e.,~$\rmax(t)\sim r_0L(t)$.
            \item The solution undergoes a quasi self-similar collapse with the
                asymptotic profile
                \begin{equation}    \label{eq:crit_ring_QSS-2}
                    \psi_{\BQ}(t,r) =
                        \frac1{L^{d/2}(t)}
                        \BQ\left( \frac {r-r_0L}L \right)
                        e^{i\int_{s=0}^{t}\frac{1}{L^{4}(s)}ds} .
                \end{equation}
            \item The blowup rate is exactly a quartic root, i.e.,
                    \begin{equation}    \label{eq:rate_14}
                        L(t)\sim\kappa\sqrt[4]{\TCrit-t},
                    \qquad \kappa>0.
                \end{equation}
%           \item The self-similar profile~$\BQ(\rho)$ is the solution of
%               \begin{equation}    \label{eq:crit_ring_ODE}
%                   \begin{gathered}
%                       -\BQ(\rho)
%                       + i\frac{\kappa^4}4 \left(
%                           \frac d2 B +(r_0+\rho) B^\prime
%                       \right)
%                   - \widetilde{\Delta_\rho^2}B + \abs{\BQ}^{8/d}B = 0
%                   \end{gathered}
%               \end{equation}
%               where~$\widetilde{\Delta_\rho^2}$ is given
%               by~\eqref{eq:crit_mangled_bi_Laplacian}.
        \end{enumerate}
    \end{subequations}
\end{conj}

\subsection{\label{ssec:critical_simulations}Simulations}

\begin{figure}
    \centering
    \subfloat{%
        \label{fig:crit_ring_rmax}%
        \includegraphics[angle=-90,clip,width=\myfw\textwidth]%
            {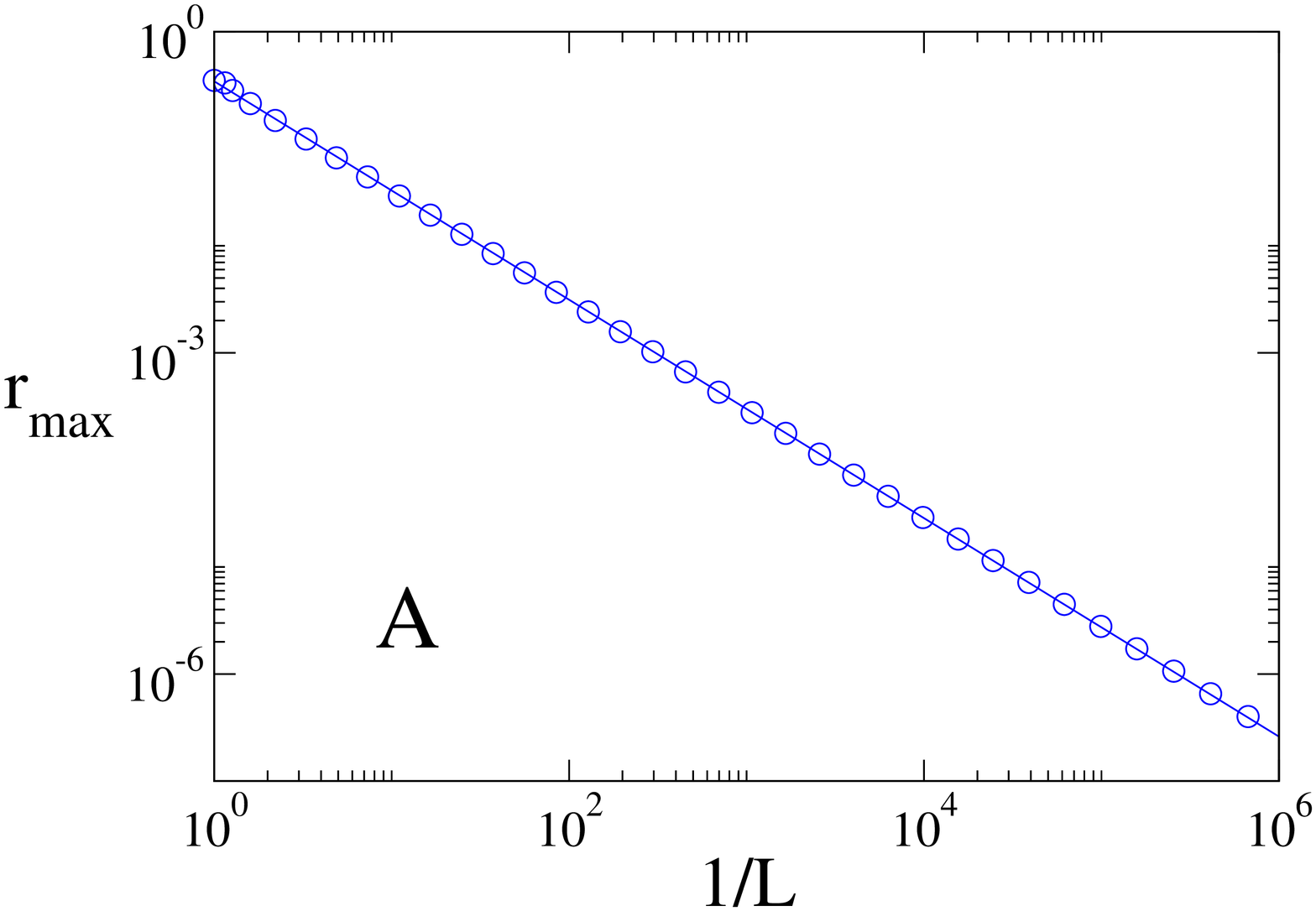}%
    }
    \subfloat{\label{fig:crit_ring_rescaled}%
        \includegraphics[angle=-90,clip,width=\myfw\textwidth]%
            {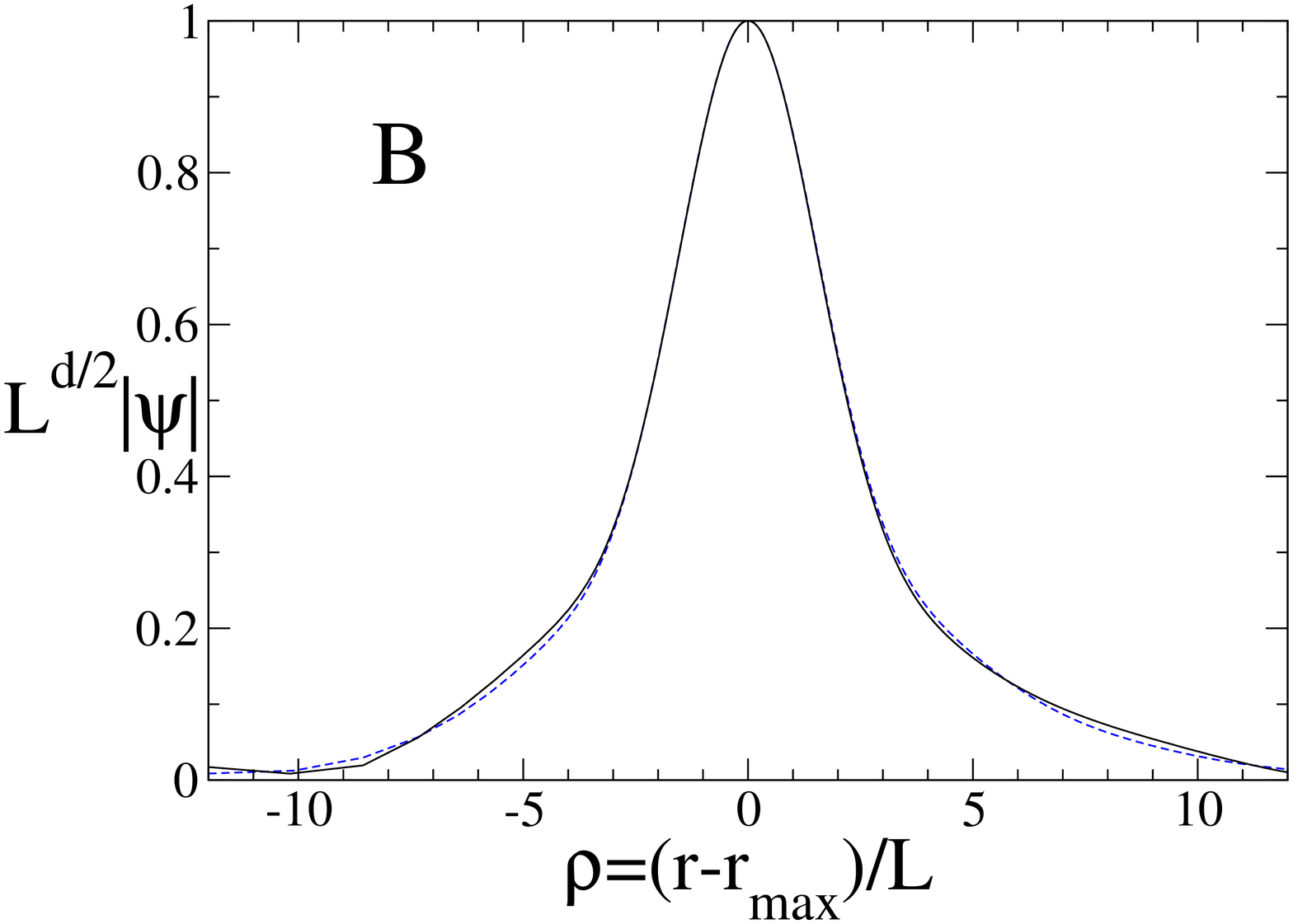}%
    }

    \mycaption{\label{fig:crit_ring_rmax_rescaled}
    Ring-type singular solution of the critical BNLS~\eqref{eq:radial_CBNLS}
    with~$d=2$.
    A)~$\rmax$ as a function of~$1/L$.
    Solid line is~$\left.\rmax=79.5 L^{1.02}\right.$.
    B) The rescaled solution, see~\eqref{eq:psi_rescaled},
    at~$L(t)=10^{-3}$ (blue dashed line) and~$L(t)=10^{-6}$ (black solid line).
    }
\end{figure}
\begin{figure}
    \centering
    \subfloat{\label{fig:crit_ring_L}%
        \includegraphics[angle=-90,clip,width=\myfw\textwidth]%
            {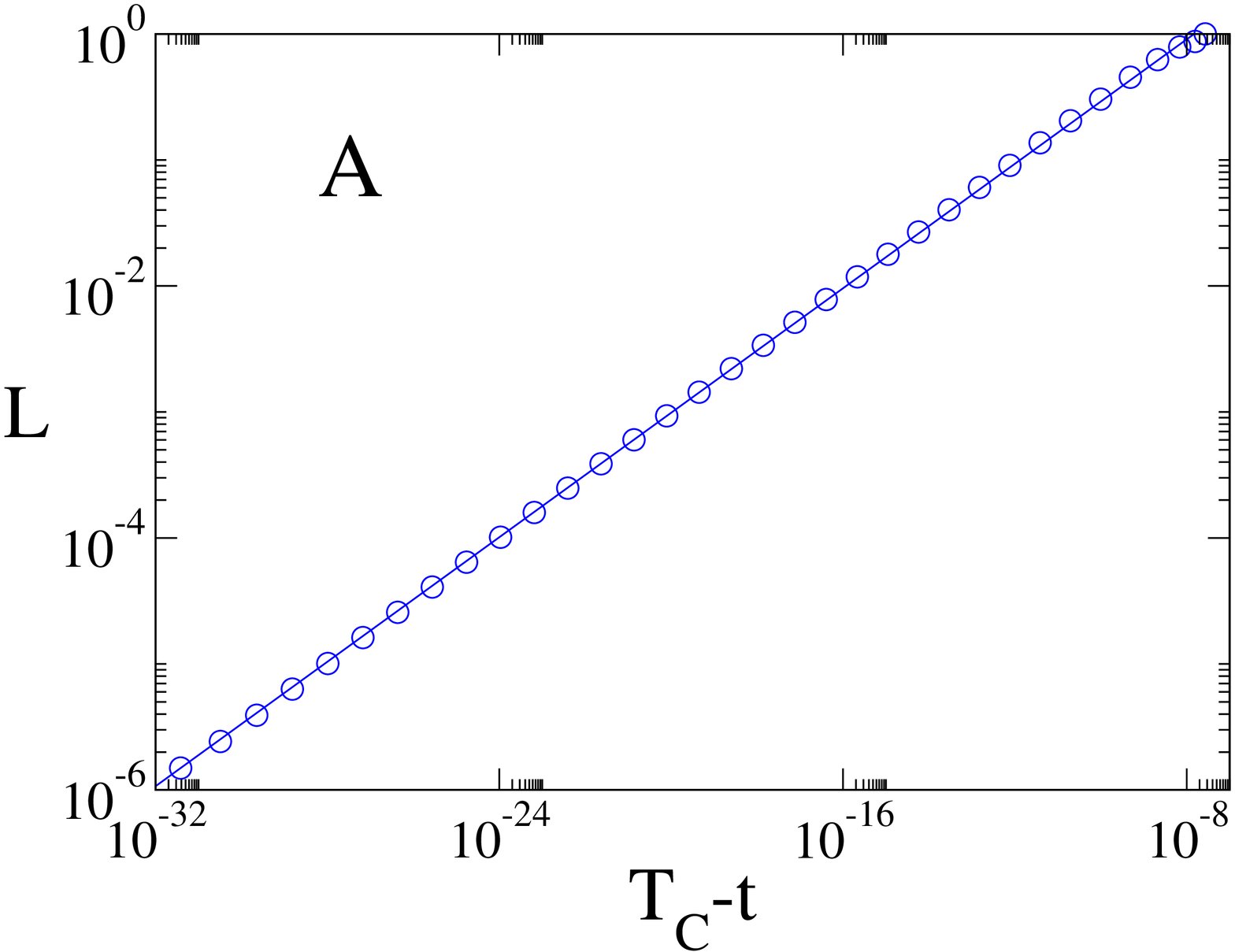}%
    }
    \subfloat{\label{fig:crit_ring_LtL3}%
        \includegraphics[angle=-90,clip,width=\myfw\textwidth]%
            {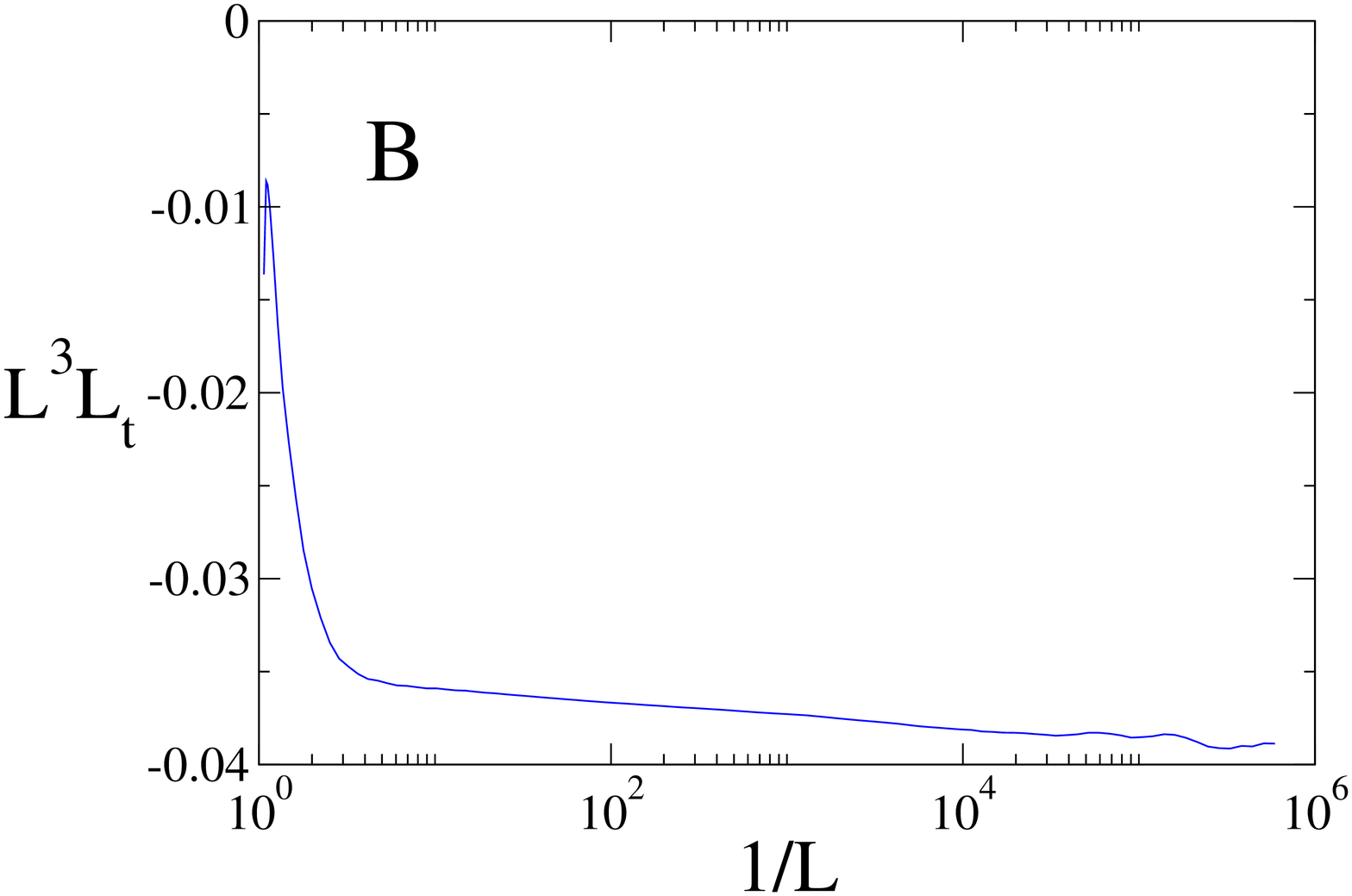}%
    }

    \mycaption{\label{fig:crit_ring_L_LtL3}
    Blowup rate of the solution of Figure~\ref{fig:crit_ring_rmax_rescaled}.
    A) Solution width~$L$ as a function of~$\left(\TCrit-t\right)$ on a
    logarithmic scale (circles).
    Solid line is~$\left.L=0.433(\TCrit-t)^{0.2476}\right.$.
    B) $L_tL^3$ as a function of~$1/L$.
    }
\end{figure}
The critical BNLS~\eqref{eq:radial_CBNLS} with~$d=2$ was solved with the
initial condition~$\psi_0=2.5 e^{(r-10)^2}$.
The simulation was run up to a focusing level of~$L(t)=10^{-6}$.
We next test Conjecture~\ref{conj:crit_ring_rate_profile} numerically, clause by
clause.
\begin{enumerate}
    \item Figure~\ref{fig:crit_ring_rmax} shows that the ring shrinks at a rate
        of~$\rmax(t)\sim cL^\alpha(t)$ with $\alpha\approx1.02$, which is close to
        the predicted value of~$\alpha=1$.
        % c=79.5
    \item In Figure~\ref{fig:crit_ring_rescaled} we plot the solution, rescaled
        according to~\eqref{eq:psi_rescaled}, at the focusing levels~$1/L=10^3$
        and~$1/L=10^6$.
        The two curves are indistinguishable, showing that the solution
        undergoes a quasi self-similar collapse with the~$\psi_{\BQ}$
        profile~\eqref{eq:crit_ring_QSS-2}.
    \item Figure~\ref{fig:crit_ring_L} shows that~$
            L(t)\sim 0.433(\TCrit-t)^{0.2477}.
        $
        Therefore, the calculated blowup rate is close to a quartic root.
    \item By Conjecture~\ref{conj:crit_ring_rate_profile}, the blowup rate
        of~$L(t)$ should be exactly~$1/4$, hence~$L^3L(t)$ should converge to a
        negative constant.
        However, in Figure~\ref{fig:crit_ring_LtL3}, $L^3L_t$ does not converge
        to a constant, but rather slowly decreases away from zero.
        This indicates that the blowup rate is slower than a quartic-root,
        which is in contradiction with Theorem~\ref{thrm:low-bound}.
        There are two possible explanations for this:
        \begin{enumerate}
            \item It may be that the numerical finding that~$\alpha$ is
                slightly above~$1$ and~$p$ is slightly below~$1/4$ is an
                artifact of our numerical method, see Section~\ref{sec:SGR}.
                Indeed, in all our simulations for~$1/2\le\alpha<1$ in
                Figures~\ref{fig:shrinking_rings_alpha}
                and~\ref{fig:shrinking_rings_p}, the calculated values of the
                shrinking rate~$\alpha$ were all slightly above~$\alphab$, and
                the blowup rates were slightly below~$\frac{1}{3+\alphab}$.
                In those cases, however, these small differences did not change
                the qualitative behavior of the solution.
                In contrast, a small increment (whether numerical or genuine)
                from~$\alpha=1$ will drastically change the dynamics, from
                equal-rate ring-type solutions into peak-type solutions.
            \item It may be that ring-type solutions of the critical BNLS are
                only meta-stable, having shrinking rates~$\alpha>1$ and
                a blowup rate slower than~$1/4$.
                This does not contradict with Theorem~\ref{thrm:low-bound},
                since in this case the ring-type solutions will eventually
                evolve into peak-type solutions with a different blowup rate.
        \end{enumerate}
        We do not know which of the two options is true.
\end{enumerate}

\section{\label{sec:SGR}Numerical Method: Adaptive mesh construction}

In this study, we computed singular solutions of the BNLS
equation~\eqref{eq:radial_BNLS}.
These solutions become highly-localized, so that the spatial scale-difference
between the singular region~$r-\rmax = {\cal O}(L)$ and the exterior regions
can be as large as~${\cal O}(1/L)\sim 10^{10}$.
In order to resolve the solution at both the singular and non-singular regions,
we use an adaptive grid.

We generate the adaptive grids using the {\em Static Grid Redistribution}~(SGR)
method, which was first introduced by Ren and Wang~\cite{Ren-00}, and later
simplified and improved by Gavish and Ditkowsky~\cite{SGR-08}.
Using this approach, the solution is allowed to propagate (self-focus) until it
becomes under-resolved.
At this stage, a new grid, with the same number of grid-points, is generated
using De'Boors `equidistribution principle', wherein the grid points~$\{r_m\}$
are spaced such that a certain weight function~$w_1[\psi]$ is equidistributed,
i.e., that \[
	\int_{r=r_m}^{r_{m+1}} w_1\left[ \psi(r) \right]dr
	= \text{const},
\]
see~\cite{Ren-00,SGR-08} for details.

\begin{algorithm}[H]
    \begin{center}
      \begin{enumerate}
          \item Find a nonlinear coordinate
              transformation~$r(x):[0,1]\to[0,R]$, under which the weight
              function~$w\left[\psi(r(x))\right]$ becomes uniformly
              distributed.
          \item Transform the solution and equation to the new coordinate
              system.
              For example, the second spatial derivative in the NLS transforms
              as \[
                  \psi_{rr} \mapsto \psi_{xx} r_x^2 + \psi_x r_{xx}.
              \]
          \item Approximate the equation on a uniform 
              grid~$\left\{ x_m \right\}$, using standard finite-differences
              (or another method of choice).
      \end{enumerate}
      \end{center}
    \caption{
        The SGR method, as implemented in~\cite{SGR-08}.
    }
    \label{alg:old-SGR}
\end{algorithm}
The method of~\cite{SGR-08} is given in Algorithm~\ref{alg:old-SGR}.
Note that, since~$r(x)$ is nonlinear, the mapping of the derivatives
of~$\psi$ (step 2) involves nonlinear combinations of the derivatives of~$r$.
This is not a great problem for the NLS, which has only second-order
derivatives, but becomes much messier for the biharmonic
operator~\eqref{eq:radial_bi_Laplacian}, with its many high-order derivatives
of~$\psi$.

Therefore, in this study we implement a simplified version of the method
of~\cite{SGR-08}, which is given in Algorithm~\ref{alg:new_SGR}, which uses a
non-uniform grid in the old-coordinate system, and thereby dispenses with the
need for transforming {\em the equation}, and is much easier to implement in
the biharmonic case.
\begin{algorithm}[H]
    \begin{center}
        \begin{enumerate}
            \item Find a nonlinear coordinate
                transformation~$r(x):[0,1]\to[0,R]$, under which the weight
                function~$w\left[\psi(r(x))\right]$ becomes uniformly
                distributed.
            \item Create the uniform grid in the transformed 
                system~$\left\{x_m \right\}$.
            \item Create the (highly) non-uniform grid~$r_m=r(x_m)$ in the
                original (physical) coordinate system.
            \item On the non-uniform grid, approximate the equation using
                standard (non-uniform) finite-differences.
      \end{enumerate}
    \end{center}
    \caption{\label{alg:new_SGR}%
        The SGR method, as implemented in this work.
    }
\end{algorithm}
\noindent
We use a third-order accurate finite-difference
approximation of the radial biharmonic
operator~\eqref{eq:radial_bi_Laplacian}, with a seven-point stencil.

\begin{figure}
	\centering
	\subfloat[old method,~$L=10^{-6}$]{\label{fig:w3_effect_old}%
		\includegraphics[clip,width=\myfw\textwidth]%
			{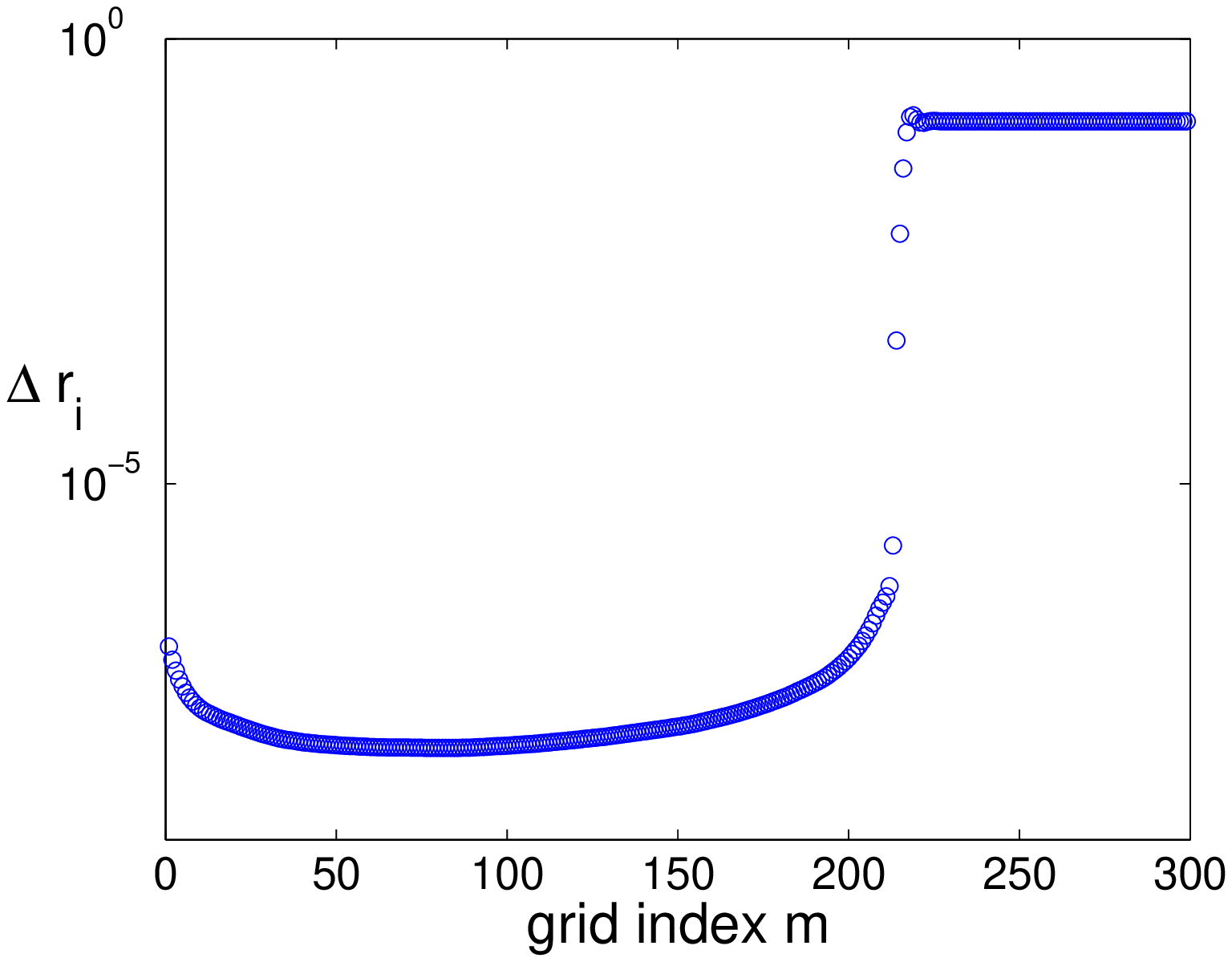}%
	}
	\subfloat[new method,~$L=10^{-12}$]{\label{fig:w3_effect_new}%
		\includegraphics[clip,width=\myfw\textwidth]%
			{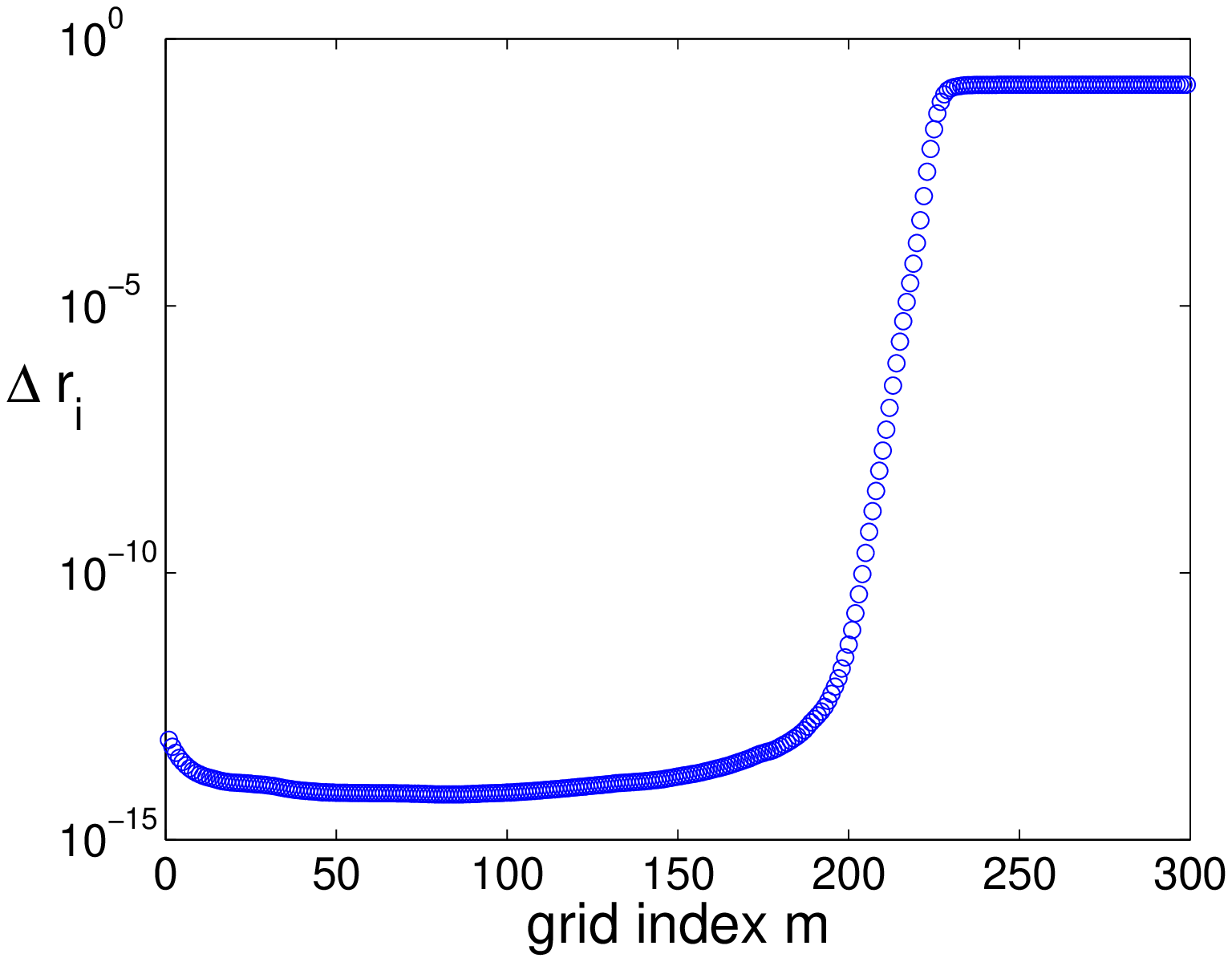}%
	}

	\mycaption{\label{fig:w3_effect}%
		The grid-spacing~$\Delta r_m$ obtained using the SGR method
		of~\cite{SGR-08} for a peak-type singular solution of the BNLS.
		A) The grid generated the original method of~\cite{SGR-08} at focusing
		level of~$L=10^{-6}$. 
		The Singular and non-singular regions are well-resolved, but the
		transition region~$\Delta r_m$ displays a discontinuity.
		At this point, the finite difference operator becomes ill-conditioned.
		B) same as (A), after adding the new penalty function~$w_3$, at focusing
		level~$L=10^{-12}$.
		Even at this much larger focusing level, the transition region is now
		well resolved.
	}
	%}
\end{figure}
The method in~\cite{SGR-08} allows control of the fraction of grid points
that migrate into the singular region, preventing under-resolution at the
exterior regions.
This is done by using a weight-function~$w_2$, which penalizes large inter-grid
distances.
However, we found that this numerical mechanism, while necessary, is
insufficient for our purposes.
In order to understand the reason, let us consider the grid-point spacings~$
	\Delta r_m = r_{m+1}-r_m
$.
Using the method of~\cite{SGR-08} with both~$w_1$ and~$w_2$ causes a very sharp 
bi-partition of the grid points -- to those inside the singular region, whose
spacing is determined by~$w_1$ and is~$\Delta r_m=\mathcal{O}(L)$,
and to those outside the singular region, whose spacing is determined by~$w_2$ 
and is~$\Delta r_m=\mathcal{O}(1)$, see Figure~\ref{fig:w3_effect_old}.
Inside each of these regions, the finite difference approximation we use is well
conditioned.
However, at the transition between these two regions, the finite-difference stencil,
seven-points in width, spans grid-spacings with~$\mathcal{O}(1/L)$
scale-difference --- leading to under-resolution which completely violates the
validity of the finite-difference approximation.

In order to overcome this limitation, we improve the algorithm of~\cite{SGR-08}
by adding a third weight function \[
w_3 (r_m) = \sqrt{	1+ \frac{\abs{\Delta^2r_m}}{\Delta r_m} },
\]
which penalizes the second-difference~$\Delta^2 r_m = \Delta r_{m+1}-\Delta r_m$
operator of the grid locations, allowing for a smooth transition between the
singular region and the non-singular region, see Fig~\ref{fig:w3_effect_new}.

On the sequence of grids, the equations are solved using a
Predictor-Corrector Crank-Nicholson scheme, which is second-order in time.

\subsection*{Acknowledgments} 
We thank Nir Gavish for useful discussions.
This research was partially supported by grant \#123/2008 from the Israel
Science Foundation (ISF).

%alpha  apalike  plain.bst   unsrt.bst

%\bibliographystyle{unsrt}

%\bibliographystyle{alpha}
%\bibliography{NLS,BNLS}
\newcommand{\etalchar}[1]{$^{#1}$}

\end{document}